\documentclass[preprint,3p,12pt]{elsarticle}
\usepackage{amsmath}
\usepackage{stmaryrd}
\usepackage{bbding}
\usepackage{dcolumn}
\usepackage{graphicx}
\usepackage{amsfonts}
\usepackage{amssymb}
\usepackage{psfrag}
\usepackage{wrapfig}
\usepackage{subfigure}
\usepackage{makeidx}
\usepackage{bm}
\usepackage{epsf}
\usepackage{epsfig}
\usepackage{setspace}
\usepackage{color}

\newtheorem{rmk}{Remark}

\begin{document}
\title{High-order gas-kinetic scheme with three-dimensional WENO reconstruction for the Euler and Navier-Stokes solutions}
\author[BNU]{Liang Pan\corref{cor}}
\ead{panliang@bnu.edu.cn}
\author[HKUST1,HKUST2,HKUST3]{Kun Xu}
\ead{makxu@ust.hk}
\address[BNU]{School of Mathematical Sciences, Beijing Normal University, Beijing, China}
\address[HKUST1]{Department of mathematics, Hong Kong University of Science and Technology, Kowloon, Hong Kong}
\address[HKUST2]{Department of Mechanical and Aerospace Engineering, Hong Kong University of Science and Technology, Kowloon, Hong Kong}
\address[HKUST3]{Shenzhen Research Institute, Hong Kong University of Science and Technology, Shenzhen,
China} \cortext[cor]{Corresponding author}

\begin{abstract}
In this paper, a simple and efficient third-order weighted
essentially non-oscillatory (WENO) reconstruction is developed for
three-dimensional flows, in which the idea of two-dimensional
WENO-AO scheme on unstructured meshes \cite{WENO-ao-3} is adopted.
In the classical finite volume type WENO schemes, the linear weights
for the candidate stencils are obtained by solving linear systems at
Gaussian quadrature points of cell interface. For the
three-dimensional scheme, such operations at twenty-four Gaussian
quadrature points of a hexahedron would reduce the efficiency
greatly, especially for the moving-mesh computation. Another
drawback of classical WENO schemes is the appearance of negative
weights with irregular local topology, which affect the robustness
of spatial reconstruction. In such three-dimensional WENO-AO scheme,
a simple strategy of selecting big stencil and sub-stencils for
reconstruction is proposed. With the reconstructed quadratic
polynomial from big stencil and linear polynomials from
sub-stencils, the linear weights are chosen as positive  numbers
with the requirement that their sum equals one and be independent of
local mesh topology. With such WENO reconstruction, a high-order
gas-kinetic scheme (HGKS) is developed for both three-dimensional
inviscid and viscous flows. Taken the grid velocity into account,
the scheme is extended into the moving-mesh computation as well.
Numerical results are provided to illustrate the good performance of
such new finite volume WENO schemes. In the future, such WENO
reconstruction will be extended to the unstructured meshes.
\end{abstract}

\maketitle

\section{Introduction}
In recent decades, there have been continuous interests and efforts
on the development of high-order schemes for compressible flows.
There have been a gigantic number of publications about the
introduction and survey of high-order schemes, including
discontinuous Galerkin (DG) \cite{DG2,DG3}, essentially
nonoscillatory (ENO), weighted ENO (WENO) schemes, etc. In this
paper, we mainly focus on the finite volume type schemes. The ENO
schemes were proposed in \cite{ENO-1,ENO-2} and successfully applied
to solve the hyperbolic conservation laws and other convection
dominated problems. The ``smoothest" stencil is selected among
several candidates to achieve high-order accuracy in the smooth
region and keep essentially non-oscillatory near discontinuities.
For unstructured meshes, the ENO scheme was developed as well
\cite{ENO-3}. Following the ENO scheme, the WENO schemes
\cite{WENO-Liu,WENO-JS,WENO-M,WENO-Z} were developed. With the
nonlinear convex combination of candidate polynomials, WENO scheme
achieves higher order of accuracy and keeps non-oscillatory property
essentially. Compared with ENO scheme, the WENO schemes improve
robustness, smoothness of fluxes, steady-state convergence and
efficiency in the computation.  On the unstructured meshes, the WENO
schemes were also developed \cite{un-WENO1}. Similar with the
one-dimensional WENO scheme, the high-order of accuracy is obtained
by the combination of lower order polynomials. However, its
successful application is limited by the appearance of negative
linear weights and very large linear weights, which appears commonly
on the unstructured meshes. For a mesh that is close to the regular
meshes, such WENO scheme works well by a regrouping approach to
avoid negative weights. However, for a mesh with lower quality, the
large linear weights appears and WENO schemes become unstable even
for the smooth flows. In order to avoid the negative linear weights
and very large linear weights, many WENO schemes were proposed
\cite{un-WENO2,un-WENO3,un-WENO4}. Instead of concentrating on the
reconstruction of interface values, there exist another class of
WENO methods to reconstruct a polynomial inside each cell based on
all stencils, which is also named as the WENO with adaptive order
(WENO-AO) method \cite{WENO-ao-1,WENO-ao-2,WENO-ao-3}.  The linear
weights are artificially set to be positive numbers with the
requirement that their sum equals to one. With the non-linear weights,
the WENO-AO schemes could achieve the optimal order of accuracy in
smooth region, and automatically approach to the smoothest quadratic
sub-stencil in discontinuous region using the same stencils from
original WENO scheme. The independence of linear weights on local
topology not only improve the efficiency, but also reduces the
complexity of the classical WENO scheme.

In the past decades, the gas-kinetic schemes (GKS) based on the
Bhatnagar-Gross-Krook (BGK) model \cite{BGK-1,BGK-2} have been
developed systematically for the computations from low speed flow to
supersonic one \cite{GKS-Xu1,GKS-Xu2}. Different from the numerical
methods based on Riemann fluxes \cite{Riemann-appro}, a
time-dependent gas distribution function is provided at the cell
interface for inviscid and viscous terms together. With such spatial and temporal
coupled gas distribution function, the one-stage third-order GKS was
developed \cite{GKS-high-Li}, which integrates the flux function
over a time step analytically without employing the multi-stage
Runge-Kutta time stepping techniques \cite{TVD-RK}. However, with
the one-stage gas evolution model, the formulation of GKS becomes
very complicated for the further improvement \cite{GKS-high-Niu}.
Recently, based on the time-dependent flux function, a two-stage
fourth-order method was developed for Lax-Wendroff type flow solvers
\cite{GRP1,GKS-Xu2}, particularly for the hyperbolic conservation
laws \cite{GRP-high-1,GRP-high-2,GKS-high-1,GKS-high-2}. With the
two-stage temporal discretization, a reliable framework was provided
for developing fourth-order and even higher-order accuracy \cite{zhao-gks8}. For the
construction of high-order scheme, a spatial-temporal coupled
evolution model becomes important, and the delicate flow structures
depend on the quality of flow solvers \cite{GKS-high-3,GKS-high-2}.
With the dimensional-by-dimensional WENO reconstruction, the
high-order gas-kinetic scheme has been extended to three-dimensional
computation with the structured meshes \cite{GKS-high-5}, especially
for the direct numerical simulation for the compressible isotropic
turbulence \cite{GKS-high-4}.

To simulate the flow with complicated geometry, a high-order
gas-kinetic scheme was proposed with the unstructured WENO
reconstruction, and extended to the moving-mesh computation
\cite{GKS-ALE}. The accuracy and geometric conservation law are well
preserved even with the largely deforming mesh. However, choosing
sub-stencils from big stencil and solving linear weights at Gaussian
quadrature points would make the reconstruction extremely
complicated for three-dimensional flows. In this paper, a simple and
efficient third-order WENO-AO scheme is developed for
three-dimensional flows to overcome the drawbacks above, in which
the three-dimensional structured mesh is considered for simplicity.
With the such reconstruction, a high-order gas-kinetic scheme (HGKS)
is developed for both Euler and Navier-Stokes solutions. In such
three-dimensional WENO-AO scheme, the strategy of selecting big
stencil and candidates of  sub-stencils for reconstruction is proposed.
Based on the reconstructed quadratic polynomial for big stencil and
linear polynomials for sub-stencils, the spatial independent linear
weights are used, which have fixed values and become positive. With
the smooth indicator, the nonlinear weights can be constructed.
Meanwhile, the point-values and slopes for non-equilibrium part of
gas distribution function can be reconstructed at all Gaussian
quadrature points. Through particle colliding procedure, the
point-values and slopes for equilibrium part are obtained
simultaneously and an extra reconstruction for equilibrium state in
the classical HGKS is avoided. Taken the grid velocity into account,
such scheme can be also extended into the moving-mesh computation.
For the mesh with non-coplanar vertexes, which is commonly generated
in the moving-mesh computation, the trilinear interpolation is used
to parameterize the hexahedron, and the bilinear interpolation is
used to parameterize the interface of hexahedron. Extensive
numerical results are provided to illustrate the good performance of
such new finite volume WENO schemes. The optimal order of accuracy
in smooth regions can be obtained, and the strong discontinuities
are also well captured.

This paper is organized as follows. In Section 2, third-order WENO
reconstruction is introduced. The high-order gas-kinetic scheme is
presented in Section 3. In Section 4, we present the extension to
ALE framework. Section 5 includes numerical examples to validate the
current algorithm. The last section is the conclusion.

\section{Third-order WENO reconstruction}
\begin{figure}[!h]
\centering
\includegraphics[width=0.5\textwidth]{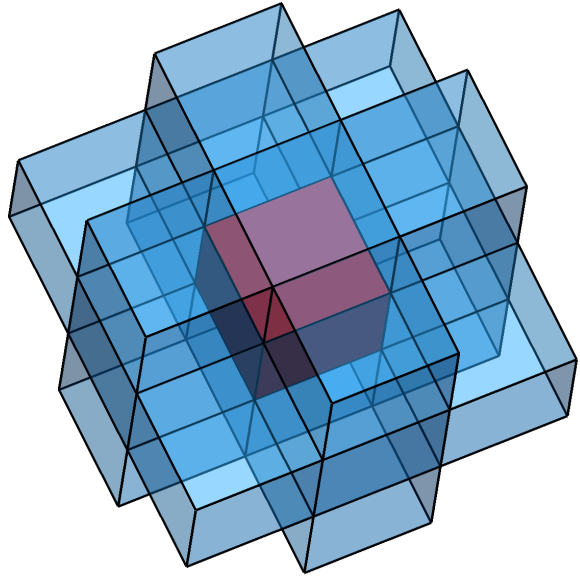}
\caption{\label{schematic} Stencils of cell $\Omega_{ijk}$ (red
cubic) for three-dimensional reconstruction.}
\end{figure}

In this section, an efficient and simple WENO scheme is proposed for
three-dimensional flows, and the idea comes from the
two-dimensional WENO-AO scheme on unstructured meshes \cite{WENO-ao-3}.
For simplicity, the reconstruction is developed on structured meshes
and the extension to unstructured meshes will be developed in the
future. For a piecewise smooth function $Q(\boldsymbol{x})$ over
cell $\Omega_{ijk}$, a polynomial $P_0(\boldsymbol{x})$ with degree
$r$ can be constructed to approximate $Q(\boldsymbol{x})$ as follows
\begin{equation*}
\displaystyle
p(\boldsymbol{x})=Q(\boldsymbol{x})+\mathcal{O}(h^{r+1}),
\end{equation*}
where $|\Omega_{ijk}|$
is volume of $\Omega_{ijk}$ and $h\sim|\Omega_{ijk}|^{1/3}$ is the cell size. In order to achieve third-order
accuracy and satisfy conservative property, the following quadratic
polynomial is constructed
\begin{equation}\label{qua-def}
P_0(\boldsymbol{x})=Q_{ijk}+\sum_{|n|=1}^2a_np_{n}(\boldsymbol{x}),
\end{equation}
where $Q_{ijk}$ is the cell averaged variables of $Q(\boldsymbol{x})$ over $\Omega_{ijk}$,
$n=(n_1, n_2, n_3)$, $|n|=n_1+n_2+n_3$ and
\begin{align*}
\displaystyle
p_{n}(\boldsymbol{x})=x^{n_1}y^{n_2}z^{n_3}-\frac{1}{\left|
\Omega_{ijk}
\right|}\iiint_{\Omega_{ijk}}x^{n_1}y^{n_2}z^{n_3}\text{d}V.
\end{align*}
In order to fully determine this polynomial, a big stencil $S$ for
$\Omega_{ijk}$, which is shown in Fig.\ref{schematic}, is selected
as follows
\begin{align*}
S=\{\Omega_{i+i_0,j+j_0,k+k_0},~i_0,j_0,k_0=-1,0,1,~i_0\cdot
j_0\cdot k_0\neq\pm1\}.
\end{align*}
The following constrains need to be satisfied  in the
big stencil
\begin{align*}
\frac{1}{\left|\Omega_{i'j'k'}\right|}\iiint_{\Omega_{i'j'k'}}P_0(\boldsymbol{x})\text{d}V=Q_{i'j'k'},~\Omega_{i'j'k'}\in
S.
\end{align*}
An over determined linear system will be generated, and the coefficients
$a_n$ in Eq.\eqref{qua-def} can be determined by the least square
method.

Twenty-four candidate sub-stencils are selected from the big stencil as
well, and one third of them are given as follows
\begin{align*}
S_1=\{\Omega_{i,j,k},\Omega_{i+1,j,k},\Omega_{i,j+1,k},\Omega_{i+1,j+1,k},\Omega_{i,j,k+1}\},\\
S_2=\{\Omega_{i,j,k},\Omega_{i-1,j,k},\Omega_{i,j+1,k},\Omega_{i-1,j+1,k},\Omega_{i,j,k+1}\},\\
S_3=\{\Omega_{i,j,k},\Omega_{i+1,j,k},\Omega_{i,j-1,k},\Omega_{i+1,j-1,k},\Omega_{i,j,k+1}\},\\
S_4=\{\Omega_{i,j,k},\Omega_{i-1,j,k},\Omega_{i,j+1,k},\Omega_{i-1,j-1,k},\Omega_{i,j,k+1}\},\\
S_5=\{\Omega_{i,j,k},\Omega_{i+1,j,k},\Omega_{i,j+1,k},\Omega_{i+1,j+1,k},\Omega_{i,j,k-1}\},\\
S_6=\{\Omega_{i,j,k},\Omega_{i-1,j,k},\Omega_{i,j+1,k},\Omega_{i-1,j+1,k},\Omega_{i,j,k-1}\},\\
S_7=\{\Omega_{i,j,k},\Omega_{i+1,j,k},\Omega_{i,j-1,k},\Omega_{i+1,j-1,k},\Omega_{i,j,k-1}\},\\
S_8=\{\Omega_{i,j,k},\Omega_{i-1,j,k},\Omega_{i,j+1,k},\Omega_{i-1,j-1,k},\Omega_{i,j,k-1}\}.
\end{align*}
Symmetrically, another two thirds of the sub-stencils can be fully
given. For each sub-stencil, a linear polynomial
$P_m(\boldsymbol{x}), m=1,...,24$ is constructed
\begin{equation}\label{linear-def}
P_m(\boldsymbol{x})=Q_{ijk}+\sum_{|n|=1}b_n^mp_{n}(\boldsymbol{x}).
\end{equation}
The following constrains need to be satisfied for cells in the
candidate sub-stencils
\begin{align*}
\frac{1}{\left|\Omega_{i'j'k'}\right|}\iiint_{\Omega_{i'j'k'}}P_m(\boldsymbol{x})\text{d}x\text{d}y=Q_{i'j'k'},
~\Omega_{i'j'k'}\in S_m.
\end{align*}
The coefficients $b_n^m$ in Eq.\eqref{linear-def} can be determined
by the least square method as well.

With the reconstructed polynomials, the quadratic polynomial
$P_0(\boldsymbol{x})$ can be rearranged as follows
\begin{equation*}
P_0(\boldsymbol{x})=\gamma_0(\frac{1}{\gamma_0}P_0(\boldsymbol{x})-\sum_{m=1}^{24}\frac{\gamma_m}{\gamma_0}P_m(\boldsymbol{x}))+\sum_{m=1}^{24}\gamma_mP_m(\boldsymbol{x}),
\end{equation*}
where $\gamma_m>0, m=0,...,24$ are the linear weights,
$\gamma_0=1-24\Gamma, \gamma_i=\Gamma, i=1,...,24$. In the
computation, $\Gamma=0.0025$ is suggested. To deal with the
discontinuities, the non-linear weights are introduced. The
reconstructed point-value to approximate $Q(\boldsymbol{x}_{G})$ at
Gaussian quadrature point is given as
\begin{equation*}
\widetilde{P}_0(\boldsymbol{x}_{G})=\overline{\omega}_0(\frac{1}{\gamma_0}P_0(\boldsymbol{x}_{G})-\sum_{m=1}^{24}\frac{\gamma_m}{\gamma_0}P_m(\boldsymbol{x}_{G}))+\sum_{m=1}^{24}\overline{\omega}_mP_m(\boldsymbol{x}_{G}).
\end{equation*}
The non-linear weights $\omega_m$ and normalized non-linear weights
$\overline{\omega}_m$ are defined as
\begin{align*}
\overline{\omega}_{m}=\frac{\omega_{m}}{\sum_{m=0}^{24} \omega_{m}},
m=0,...,24
\end{align*}
and
\begin{align}\label{nonlinear}
\omega_{m}=\gamma_{m}\Big[1+\Big(\frac{\tau}{\beta_{m}+\epsilon}\Big)\Big],
\end{align}
where $\epsilon$ is a small positive number. The smooth indicator
$\beta_{m}$ is given as
\begin{align}\label{indicator}
\beta_m=\sum_{|l|=1}^{r_m}|\Omega_{ijk}|^{\frac{2|l|}{3}-1}\int_{\Omega_{ijk}}\Big(\frac{\partial^lP_m}{\partial_x^{l_1}\partial_y^{l_2}\partial_z^{l_3}}(x,y,z)\Big)^2\text{d}V,
\end{align}
where $r_0=2$ and $r_m=1, m\neq0$. In order to achieve the optimal
order of accuracy, the parameter $\tau$ is chosen as
\begin{align}\label{para}
\tau=\sum_{m=1}^{24}\Big(\frac{|\beta_0-\beta_m|}{24}\Big).
\end{align}
According to the Taylor expansion, the smooth indicator
Eq.\eqref{indicator} can be rewritten as
\begin{align*}
\beta_0=&\Big(\big(\frac{\partial Q}{\partial
x}\big)^2+\big(\frac{\partial Q}{\partial
y}\big)^2+\big(\frac{\partial Q}{\partial
z}\big)^2\Big)\Big|_{\boldsymbol{x}_0}
|\Omega_{ijk}|^{2/3}(1+\mathcal{O}(|\Omega_{ijk}|^{2/3})),\\
\beta_m=&\Big(\big(\frac{\partial Q}{\partial
x}\big)^2+\big(\frac{\partial Q}{\partial
y}\big)^2+\big(\frac{\partial Q}{\partial
z}\big)^2\Big)\Big|_{\boldsymbol{x}_0}
|\Omega_{ijk}|^{2/3}(1+\mathcal{O}(|\Omega_{ijk}|^{1/3})),
\end{align*}
where $m=1,...,24$. According to the definition of the parameter
$\tau$ in Eq.\eqref{para}, the non-linear weights can be approximated
as
\begin{align*}
\overline{\omega}_{m}\sim\omega_{m}=\gamma_{m}(1+\mathcal{O}(h)),
\end{align*}
which ensures the optimal order of accuracy of the current scheme.


To examine the accuracy of new WENO reconstruction, it is supposed
that the following sufficient condition is satisfied for the
nonlinear weight
\begin{align}\label{approx-k}
\omega_{m}=\gamma_{m}(1+\mathcal{O}(h^k)).
\end{align}
For the quadratic polynomial $P_0(\boldsymbol{x})$, the error can be
written as
\begin{align*}
P_0(\boldsymbol{x}_{G})=Q(\boldsymbol{x}_{G})+A(\boldsymbol{x}_{G})h^{3}+\mathcal{O}(h^{4}),
\end{align*}
where $Q(\boldsymbol{x}_{G})$ is the exact solution at the Gaussian
quadrature point $\boldsymbol{x}_{G}$. For the linear polynomial
$P_m(\boldsymbol{x}_{G})$, the error can be written as
\begin{align*}
P_m(\boldsymbol{x}_{G})=Q(\boldsymbol{x}_{G})+B_{m}(\boldsymbol{x}_{G})h^{2}+\mathcal{O}(h^{3}),
\end{align*}
The reconstructed point value with non-linear weights can be written
as
\begin{align*}
\widetilde{P}_0(\boldsymbol{x}_{G})=&\frac{\overline{\omega}_0}{\gamma_0}P_0(\boldsymbol{x}_{G})+\sum_{m=1}^{24}(\overline{\omega}_m-\frac{\overline{\omega}_0\gamma_m}{\gamma_0})P_m(\boldsymbol{x}_{G}),\\
=&\frac{\overline{\omega}_0}{\gamma_0}(Q(\boldsymbol{x}_{G})+A(\boldsymbol{x}_{G})h^{3}+\mathcal{O}(h^{4}))+\sum_{m=1}^{24}(\overline{\omega}_m-\frac{\overline{\omega}_0\gamma_m}{\gamma_0})(Q(\boldsymbol{x}_{G})+B_{m}(\boldsymbol{x}_{G})h^{2}+\mathcal{O}(h^{3}))\\
=&Q(\boldsymbol{x}_{G})+\frac{\overline{\omega}_0}{\gamma_0}(A(\boldsymbol{x}_{G})h^{3}+\mathcal{O}(h^{4}))+\sum_{m=1}^{24}(\overline{\omega}_m-\overline{\omega}_0\frac{\gamma_m}{\gamma_0})(B_{m}(\boldsymbol{x}_{G})h^{2}+\mathcal{O}(h^{3}))\\
=&Q(\boldsymbol{x}_{G})+A(\boldsymbol{x}_{G})h^{3}+\sum_{m=1}^{24}B_{m}(\boldsymbol{x}_{G})h^{2+k}+\mathcal{O}(h^{3+k}+h^{4}).
\end{align*}
Thus, to achieve the third-order accuracy, $k=1$ is needed in
Eq.\eqref{approx-k}.

With the reconstructed polynomial, the spatial derivatives at
Gaussian quadrature point, which will be used in the gas-kinetic
solver, can be given as follows
\begin{align*}
\partial_x\widetilde{P}_0(\boldsymbol{x}_{G})=\overline{\omega}_0(\frac{1}{\gamma_0}\partial_x
P_0(\boldsymbol{x}_{G})-\sum_{m=1}^{24}\frac{\gamma_m}{\gamma_0}\partial_x
P_m(\boldsymbol{x}_{G}))+\sum_{m=1}^{24}\overline{\omega}_m\partial_x
P_m(\boldsymbol{x}_{G}),\\
\partial_y\widetilde{P}_0(\boldsymbol{x}_{G})=\overline{\omega}_0(\frac{1}{\gamma_0}\partial_y
P_0(\boldsymbol{x}_{G})-\sum_{m=1}^{24}\frac{\gamma_m}{\gamma_0}\partial_y
P_m(\boldsymbol{x}_{G}))+\sum_{m=1}^{24}\overline{\omega}_m\partial_y
P_m(\boldsymbol{x}_{G}),\\
\partial_z\widetilde{P}_0(\boldsymbol{x}_{G})=\overline{\omega}_0(\frac{1}{\gamma_0}\partial_z
P_0(\boldsymbol{x}_{G})-\sum_{m=1}^{24}\frac{\gamma_m}{\gamma_0}\partial_z
P_m(\boldsymbol{x}_{G}))+\sum_{m=1}^{24}\overline{\omega}_m\partial_z
P_m(\boldsymbol{x}_{G}).
\end{align*}

\begin{rmk}
In this section, we present one choice of big stencil and candidates of
sub-stencils, which work well for the numerical tests given in the
following section. However, such choice is not unique and may not be
optimal. For example, the following big stencil $S'$
\begin{align*}
S'=S\cup\{\Omega_{i\pm2,j,k},\Omega_{i,j\pm2,k},\Omega_{i,j,k\pm2}\}.
\end{align*}
and candidate sub-stencils $S'_i, i=1,..,8,$
\begin{align*}
S_1'=\{\Omega_{i,j,k},\Omega_{i+1,j,k},\Omega_{i,j+1,k},\Omega_{i,j,k+1}\},~S_5'=\{\Omega_{i,j,k},\Omega_{i+1,j,k},\Omega_{i,j+1,k},\Omega_{i,j,k-1}\},\\
S_2'=\{\Omega_{i,j,k},\Omega_{i-1,j,k},\Omega_{i,j+1,k},\Omega_{i,j,k+1}\},~S_6'=\{\Omega_{i,j,k},\Omega_{i-1,j,k},\Omega_{i,j+1,k},\Omega_{i,j,k-1}\},\\
S_3'=\{\Omega_{i,j,k},\Omega_{i+1,j,k},\Omega_{i,j-1,k},\Omega_{i,j,k+1}\},~S_7'=\{\Omega_{i,j,k},\Omega_{i+1,j,k},\Omega_{i,j-1,k},\Omega_{i,j,k-1}\},\\
S_4'=\{\Omega_{i,j,k},\Omega_{i-1,j,k},\Omega_{i,j+1,k},\Omega_{i,j,k+1}\},~S_8'=\{\Omega_{i,j,k},\Omega_{i-1,j,k},\Omega_{i,j+1,k},\Omega_{i,j,k-1}\}.
\end{align*}
can be used as substitutions respectively. These substitutions work
as well in the computation.
\end{rmk}

As mentioned in the introduction, the vertexes of hexahedron may
become non-coplanar during the moving mesh procedure, which
introduces difficulty to preserve the high-order accuracy and
geometric conservation law. In this paper, the trilinear
interpolation is used to describe the hexahedron with non-coplanar
vertex as follows
\begin{align*}
\boldsymbol{X}(\xi,\eta,\zeta)=\sum_{i=1}^8\boldsymbol{x}_i\psi_i(\xi,
\eta,\zeta),
\end{align*}
where  $(\xi,\eta,\zeta)\in[-1/2,1/2]^3$, $\boldsymbol{x}_i$ is the
vertex of a hexahedron and $\psi_i$ is the base function
\begin{align*}
\psi_1=\frac{1}{8}(1-2\xi)(1-2\eta)(1-2\zeta),~
\psi_2=\frac{1}{8}(1-2\xi)(1-2\eta)(1+2\zeta),\\
\psi_3=\frac{1}{8}(1-2\xi)(1+2\eta)(1-2\zeta),~
\psi_4=\frac{1}{8}(1-2\xi)(1+2\eta)(1+2\zeta),\\
\psi_5=\frac{1}{8}(1+2\xi)(1-2\eta)(1-2\zeta),~
\psi_6=\frac{1}{8}(1+2\xi)(1-2\eta)(1+2\zeta),\\
\psi_7=\frac{1}{8}(1+2\xi)(1+2\eta)(1-2\zeta),~
\psi_8=\frac{1}{8}(1+2\xi)(1+2\eta)(1+2\zeta).
\end{align*}
For the hexahedron with non-coplanar vertex, the triple integral
over the parameterized control volume can be given by
\begin{align*}
\iiint_{\Omega}x^{a}y^{b}z^{c}\text{d}x\text{d}y\text{d}z&=\iiint_{\Omega}x^{a}y^{b}z^{c}(\xi,\eta,\zeta)\Big|\frac{\partial(x,y,z)}{\partial(\xi,\eta,\zeta)}\Big|\text{d}\xi\text{d}\eta\text{d}\zeta.
\end{align*}
For simplicity, the Gaussian quadrature is used as well
\begin{align*}
\iiint_{\Omega}x^{a}y^{b}z^{c}\text{d}x\text{d}y\text{d}z
&=\sum_{l,m,n=1}^2\omega_{lmn}x^{a}y^{b}z^{c}(\xi_l,\eta_m,\zeta_n)\Big|\frac{\partial(x,y,z)}{\partial(\xi,\eta,\zeta)}\Big|_{(\xi_l,\eta_m,\zeta_n)}\Delta\xi\Delta\eta\Delta\zeta,
\end{align*}
where $\omega_{lmn}$ is quadrature weight and
$(\xi_l,\eta_m,\zeta_n)$ is the quadrature point. With such
quadrature rule, the reconstruction on the hexahedrons can be
conducted directly.

\section{High-order gas-kinetic scheme}
The three-dimensional BGK equation \citep{BGK-1,BGK-2} can be
written as
\begin{equation}\label{bgk}
f_t+uf_x+vf_y+wf_z=\frac{g-f}{\tau},
\end{equation}
where $\boldsymbol{u}=(u,v,w)$ is the particle velocity, $f$ is the
gas distribution function, $g$ is the three-dimensional Maxwellian
distribution and $\tau$ is the collision time. The collision term
satisfies the compatibility condition
\begin{equation}\label{compatibility}
\int \frac{g-f}{\tau}\psi \text{d}\Xi=0,
\end{equation}
where
$\displaystyle\psi=(\psi_1,...,\psi_5)^T=(1,u,v,w,\frac{1}{2}(u^2+v^2+w^2+\varsigma^2))^T$,
the internal variables
$\varsigma^2=\varsigma_1^2+...+\varsigma_K^2$,
$\text{d}\Xi=\text{d}u\text{d}v\text{d}w\text{d}\varsigma^1...\text{d}\varsigma^{K}$,
$\gamma$ is the specific heat ratio and  $K=(5-3\gamma)/(\gamma-1)$
is the degrees of freedom for three-dimensional flow.  According to
the Chapman-Enskog expansion for BGK equation, the macroscopic
governing equations can be derived \citep{GKS-Xu1,GKS-Xu2}. In the
continuum region, the BGK equation can be rearranged and the gas
distribution function can be expanded as
\begin{align*}
f=g-\tau D_{\boldsymbol{u}}g+\tau D_{\boldsymbol{u}}(\tau
D_{\boldsymbol{u}})g-\tau D_{\boldsymbol{u}}[\tau
D_{\boldsymbol{u}}(\tau D_{\boldsymbol{u}})g]+...,
\end{align*}
where $D_{\boldsymbol{u}}=\displaystyle\frac{\partial}{\partial
t}+\boldsymbol{u}\cdot \nabla$. With the zeroth-order truncation
$f=g$, the Euler equations can ba obtained. For the first-order
truncation
\begin{align*}
f=g-\tau (ug_x+vg_y+wg_z+g_t),
\end{align*}
the Navier-Stokes equations can ba obtained.

Taking moments of the BGK equation Eq.\eqref{bgk} and integrating
with respect to space, the finite volume scheme can be expressed as
\begin{align*}
\frac{\text{d}(Q_{ijk})}{\text{d}t}=\mathcal{L}(Q_{ijk}),
\end{align*}
where
the operator $\mathcal{L}$ is defined as
\begin{equation}\label{finite}
\mathcal{L}(Q_{ijk})=-\frac{1}{|\Omega_{ijk}|}\sum_{p=1}^6\iint_{\Sigma_{p}}F(Q)\cdot\boldsymbol{n}\text{d}\sigma,
\end{equation}
where $\Sigma_{p}$ is one cell interface of $\Omega_{ijk}$ and $\boldsymbol{n}$ is the outer normal direction.
A two-stage fourth-order time-accurate discretization was developed
for  Lax-Wendroff flow solvers with the generalized Riemann problem
(GRP) solver \cite{GRP-high-1} and the gas-kinetic scheme (GKS)
\cite{GKS-high-1}. Consider the following time-dependent equation
\begin{align*}
\frac{\partial Q}{\partial t}=\mathcal {L}(Q),
\end{align*}
with the initial condition at $t_n$, i.e.,
\begin{align*}
Q(t=t_n)=Q^n,
\end{align*}
where $\mathcal {L}$ is an operator for spatial derivative of flux,
the state $Q^{n+1}$ at $t_{n+1}=t_n+\Delta t$  can be updated with
the following formula
\begin{align*}
&Q^*=Q^n+\frac{1}{2}\Delta t\mathcal {L}(Q^n)+\frac{1}{8}\Delta
t^2\frac{\partial}{\partial
t}\mathcal{L}(Q^n), \\
Q^{n+1}=&Q^n+\Delta t\mathcal {L}(Q^n)+\frac{1}{6}\Delta
t^2\big(\frac{\partial}{\partial
t}\mathcal{L}(Q^n)+2\frac{\partial}{\partial
t}\mathcal{L}(Q^*)\big).
\end{align*}
It can be proved that for hyperbolic equations the above temporal
discretization provides a fourth-order time accurate solution for
$Q^{n+1}$. According to the definition of operator $\mathcal{L}$
Eq.\eqref{finite}, the numerical fluxes and its spatial derivative
is needed.

For each face of hexahedron with non-coplanar vertex, the trilinear
interpolation reduces to a bilinear interpolation. For the interface
with $\xi=1/2$, the coordinate is defined as
\begin{align}\label{bilinear}
\boldsymbol{X}(\eta,\zeta)=\sum_{i=1}^4\boldsymbol{x}_i\phi_i(\eta,\zeta),
\end{align}
where $(\eta,\zeta)\in[-1/2,1/2]^2$, $\boldsymbol{x}_i$ is the
vertex of the interface and $\phi_i$ is the base function
\begin{align*}
\phi_1=\frac{1}{4}(1-2\eta)(1-2\zeta),~
\phi_2=\frac{1}{4}(1-2\eta)(1+2\zeta),\\
\phi_3=\frac{1}{4}(1+2\eta)(1-2\zeta),~
\phi_4=\frac{1}{4}(1+2\eta)(1+2\zeta).
\end{align*}
With the parameterized cell interface, the numerical flux is
provided by the surface integral over the cell interface
\begin{align*}
\mathbb{F}_{p}(t)=\iint_{\Sigma_{p}}
F(Q)\cdot\boldsymbol{n}\text{d}\sigma=\int_{-1/2}^{1/2}\int_{-1/2}^{1/2}
F(Q(\boldsymbol{X}(\eta,\zeta)))\cdot\boldsymbol{n}\|\boldsymbol{X}_\eta\times\boldsymbol{X}_\zeta\|\text{d}\eta\text{d}\zeta.
\end{align*}
To achieve the spatial accuracy, the Gaussian quadrature is used for
the numerical flux above
\begin{align}\label{flux-x}
\mathbb{F}_{p}(t)=\sum_{m_1,m_2=1}^2\omega_{m_1,m_2}
F_{m_1,m_2}(t)\|\boldsymbol{X}_\eta\times\boldsymbol{X}_\zeta\|_{m_1,m_2}\Delta\eta\Delta\zeta,
\end{align}
where  $\omega_{m_1,m_2}$ is Gaussian quadrature weight, and
$F_{m_1,m_2}(t)$ is the numerical flux at the Gaussian quadrature
point, which can be obtained by taking moments of the gas distribution function
\begin{align}\label{flux-gau}
F_{m_1,m_2}(t)=\left(
\begin{array}{c}
F^{\rho}_{m_1,m_2}  \\
F^{\rho U}_{m_1,m_2}  \\
F^{\rho V}_{m_1,m_2}  \\
F^{\rho W}_{m_1,m_2}  \\
F^{\rho E}_{m_1,m_2}  \\
\end{array}\right)=\int\psi f(\boldsymbol{x}_{m_1,m_2},t,\boldsymbol{u},\xi)\boldsymbol{u}\cdot
\boldsymbol{n}_{m_1,m_2}\text{d}\Xi,
\end{align}
where $\boldsymbol{x}_{m_1,m_2}$ is the quadrature point and
$\boldsymbol{n}_{m_1,m_2}$ is the outer normal direction. In the
actual computation, the reconstruction is presented in a local
coordinate, which is given as follows
\begin{align*}
\boldsymbol{n}_x&=(\boldsymbol{X}_\eta\times\boldsymbol{X}_\zeta)/\|\boldsymbol{X}_\eta\times\boldsymbol{X}_\zeta\|,\\
\boldsymbol{n}_z&=\boldsymbol{X}_\zeta/\|\boldsymbol{X}_\zeta\|,\\
\boldsymbol{n}_y&=\boldsymbol{n}_z\times\boldsymbol{n}_x.
\end{align*}
With the reconstructed variables, the gas distribution function is
obtained at Gaussian quadrature point. The numerical flux can be
obtained by taking moments of it, and the component-wise form can be
written as
\begin{align}\label{flux-local}
\widetilde{F}_{m_1,m_2}(t)=\left(
\begin{array}{c}
F^{\widetilde{\rho}}_{m_1,m_2}  \\
F^{\widetilde{\rho U}}_{m_1,m_2}  \\
F^{\widetilde{\rho V}}_{m_1,m_2}  \\
F^{\widetilde{\rho W}}_{m_1,m_2}  \\
F^{\widetilde{\rho E}}_{m_1,m_2}  \\
\end{array}\right)=
\int\widetilde{u}\left(
\begin{array}{c}
1\\
\widetilde{u}  \\
\widetilde{v}  \\
\widetilde{w}  \\
\frac{1}{2}(\widetilde{u}^2+\widetilde{v}^2+\widetilde{w}^2+\xi^2) \\
\end{array}\right)
f(\boldsymbol{x}_{m_1,m_2},t,\widetilde{\boldsymbol{u}},\xi)\text{d}\widetilde{\Xi},
\end{align}
where $f(\boldsymbol{x}_{m_1,m_2},t,\widetilde{\boldsymbol{u}},\xi)$
is the gas distribution function in the local coordinate, and the particle velocity in the local coordinate is
given by
\begin{align*}
\widetilde{\boldsymbol{u}}=\boldsymbol{u}\cdot(\boldsymbol{n}_x,\boldsymbol{n}_y,
\boldsymbol{n}_z).
\end{align*}
Denote $(a_{ij})$ is the inverse of
$(\boldsymbol{n}_x,\boldsymbol{n}_y, \boldsymbol{n}_z)$, and each
component of $F_{m_1,m_2}(t)$ can be given by the combination of
fluxes in the local orthogonal coordinate
\begin{align}\label{trans}
\left\{\begin{aligned}
F^{\rho}_{m_1,m_2}=&F_{m_1,m_2}^{\widetilde{\rho}},\\
F^{\rho U}_{m_1,m_2}=&a_{11}F_{m_1,m_2}^{\widetilde{\rho U}}+a_{12}F_{m_1,m_2}^{\widetilde{\rho V}}+a_{13}F_{m_1,m_2}^{\widetilde{\rho W}},\\
F^{\rho V}_{m_1,m_2}=&a_{21}F_{m_1,m_2}^{\widetilde{\rho U}}+a_{22}F_{m_1,m_2}^{\widetilde{\rho V}}+a_{23}F_{m_1,m_2}^{\widetilde{\rho W}},\\
F^{\rho W}_{m_1,m_2}=&a_{31}F_{m_1,m_2}^{\widetilde{\rho U}}+a_{32}F_{m_1,m_2}^{\widetilde{\rho V}}+a_{33}F_{m_1,m_2}^{\widetilde{\rho W}},\\
F^{\rho E}_{m_1,m_2}=&F_{m_1,m_2}^{\widetilde{\rho E}},
\end{aligned} \right.
\end{align}

With the integral solution of BGK equation, the gas distribution
function in Eq.\eqref{flux-local} can be constructed as follows
\begin{equation*}
f(\boldsymbol{x}_{m_1,m_2},t,\boldsymbol{u},\varsigma)=\frac{1}{\tau}\int_0^t
g(\boldsymbol{x}',t',\boldsymbol{u},\varsigma)e^{-(t-t')/\tau}\text{d}t'+e^{-t/\tau}f_0(-\boldsymbol{u}t,\varsigma),
\end{equation*}
where
$\widetilde{\boldsymbol{u}}=(\widetilde{u},\widetilde{v},\widetilde{w})$
is denoted as $\boldsymbol{u}=(u,v,w)$ for simplicity in this
section,  $x_{m_1,m_2}=x'+u(t-t'),
y_{m_1,m_2}=y'+v(t-t'), z_{m_1,m_2}=z'+w(t-t')$ are the trajectory of
particles, $f_0$ is the initial gas distribution function, and $g$
is the corresponding equilibrium state.
 With the reconstruction of
macroscopic variables, the second-order gas distribution function at the cell
interface can be expressed as
\begin{align}\label{flux}
f(\boldsymbol{x}_{m_1,m_2},t,\boldsymbol{u},\varsigma)=&(1-e^{-t/\tau})g_0+((t+\tau)e^{-t/\tau}-\tau)(\overline{a}_1u+\overline{a}_2v+\overline{a}_3w)g_0\nonumber\\
+&(t-\tau+\tau e^{-t/\tau}){\bar{A}} g_0\nonumber\\
+&e^{-t/\tau}g_r[1-(\tau+t)(a_{1}^{r}u+a_{2}^{r}v+a_{3}^{r}w)-\tau A^r)]H(u)\nonumber\\
+&e^{-t/\tau}g_l[1-(\tau+t)(a_{1}^{l}u+a_{2}^{l}v+a_{3}^{l}w)-\tau
A^l)](1-H(u)),
\end{align}
where the equilibrium state $g_{0}$ and corresponding conservative
variables $Q_{0}$ and spatial derivatives in the local coordinate at
the quadrature point can be determined by the compatibility
condition Eq.\eqref{compatibility}
\begin{align*}
\int\psi g_{0}\text{d}\Xi=Q_0=\int_{u>0}\psi
g_{l}\text{d}\Xi+\int_{u<0}\psi g_{r}\text{d}\Xi,
\end{align*}
and
\begin{align*}
\frac{\partial Q_{0}}{\partial {\boldsymbol{n}_x}}=\int_{u>0}\psi
a_{1}^{l} g_{l}\text{d}\Xi+\int_{u<0}\psi a_{1}^{r}  g_{r}\text{d}\Xi,\\
\frac{\partial Q_{0}}{\partial {\boldsymbol{n}_y}}=\int_{u>0}\psi
a_{2}^{l} g_{l}\text{d}\Xi+\int_{u<0}\psi a_{2}^{r}  g_{r}\text{d}\Xi,\\
\frac{\partial Q_{0}}{\partial {\boldsymbol{n}_z}}=\int_{u>0}\psi
a_{3}^{l} g_{l}\text{d}\Xi+\int_{u<0}\psi a_{3}^{r}
g_{r}\text{d}\Xi.
\end{align*}
In the classical gas-kinetic scheme, an extra reconstruction is
needed for equilibrium state. Different from the multidimensional
scheme based on dimensional-by-dimensional reconstruction, the
selection of stencil and procedure of reconstruction introduce extra
difficulties for the genuine multidimensional scheme
\cite{GKS-high-2}. The procedure above reduces the complexity
greatly. The coefficients in Eq.\eqref{flux} can be determined by
the reconstructed directional derivatives and compatibility
condition
\begin{align*}
\displaystyle \langle a_{1}^{k}\rangle=\frac{\partial
Q_{k}}{\partial \boldsymbol{n}_x}, \langle
a_{2}^{k}\rangle=\frac{\partial Q_{k}}{\partial \boldsymbol{n}_y},
\langle a_{3}^{k}\rangle&=\frac{\partial Q_{k}}{\partial
\boldsymbol{n}_z}, \langle
a_{1}^{k}u+a_{2}^{k}v+a_{3}^{k}w+A^{k}\rangle=0,\\ \displaystyle
\langle\overline{a}_1\rangle=\frac{\partial Q_{0}}{\partial
{\boldsymbol{n}_x}}, \langle\overline{a}_2\rangle=\frac{\partial
Q_{0}}{\partial {\boldsymbol{n}_y}},
\langle\overline{a}_3\rangle&=\frac{\partial Q_{0}}{\partial
{\boldsymbol{n}_z}},
\langle\overline{a}_1u+\overline{a}_2v+\overline{a}_3w+\overline{A}\rangle=0,
\end{align*}
where $k=l,r$ and $\langle...\rangle$ are the moments of the
equilibrium $g$ and defined by
\begin{align*}
\langle...\rangle=\int g (...)\psi \text{d}\Xi.
\end{align*}
More details of the gas-kinetic scheme can be found in
\cite{GKS-Xu1}.

To implement the two-stage method, the numerical fluxes and its temporal derivative are given as follows
\begin{align*}
\mathbb{F}_{p}^n&=\sum_{m_1,m_2=1}^2\omega_{m_1,m_2}
F_{m_1,m_2}^n\|\boldsymbol{X}_\eta\times\boldsymbol{X}_\zeta\|_{m_1,m_2}\Delta\eta\Delta\zeta,\\
\partial_t\mathbb{F}_{p}^n&=\sum_{m_1,m_2=1}^2\omega_{m_1,m_2}
\partial_tF_{m_1,m_2}^n\|\boldsymbol{X}_\eta\times\boldsymbol{X}_\zeta\|_{m_1,m_2}\Delta\eta\Delta\zeta,
\end{align*}
where the coefficients $F_{m_1,m_2}^n$ and $\partial_tF_{m_1,m_2}^n$
can be given by the linear combination of
$\widetilde{F}_{m_1,m_2}^n$ and
$\partial_t\widetilde{F}_{m_1,m_2}^n$ in the local coordinate
according to Eq.\eqref{trans}. To determine these coefficients, the
time dependent numerical flux can be approximate as a linear
function
\begin{align}\label{expansion-1}
\widetilde{F}_{m_1,m_2}(t)=\widetilde{F}_{m_1,m_2}^n+ \partial_t
\widetilde{F}_{m_1,m_2}^n(t-t_n).
\end{align}
Integrating Eq.\eqref{expansion-1} over $[t_n, t_n+\Delta t/2]$ and
$[t_n, t_n+\Delta t]$, we have the following two equations
\begin{align*}
\widetilde{F}_{m_1,m_2}^n\Delta t&+\frac{1}{2}\partial_t \widetilde{F}_{m_1,m_2}^n\Delta t^2 =\widehat{\mathbb{F}}(\boldsymbol{x}_{m_1,m_2},\Delta t) , \\
\frac{1}{2}\widetilde{F}_{m_1,m_2}^n\Delta t&+\frac{1}{8}\partial_t
\widetilde{F}_{m_1,m_2}^n\Delta t^2
=\widehat{\mathbb{F}}(\boldsymbol{x}_{m_1,m_2}, \Delta t/2).
\end{align*}
where
\begin{align*}
\widehat{\mathbb{F}}(\boldsymbol{x}_{m_1,m_2},\delta)
=\int_{t_n}^{t_n+\delta}\widetilde{F}_{m_1,m_2}(t)\text{d}t=\int_{t_n}^{t_n+\delta}\int
\widetilde{u}\widetilde{\psi}
f(\boldsymbol{x}_{m_1,m_2},t,\widetilde{\boldsymbol{u}},\xi)\text{d}\Xi\text{d}t.
\end{align*}
The coefficients $\widetilde{F}_{m_1,m_2}^n$ and $\partial_t
\widetilde{F}_{m_1,m_2}^n$ can be determined by solving the linear
system.   Similarly, $F_{m_1,m_2}^{n*}$ and
$\partial_tF_{m_1,m_2}^{n*}$ at the intermediate state can be
constructed as well.

\begin{rmk}
Taken in the grid velocity into account,  the high-order gas-kinetic
scheme can be extended to the moving-mesh framework. Standing on the
moving reference, the three-dimensional BGK equation Eq.\eqref{bgk}
can be modified as
\begin{equation*}
f_t+(u-U^g)f_x+(v-V^g)f_y+(w-W^g)f_z =\frac{g-f}{\tau},
\end{equation*}
where $\boldsymbol{U}^{g}=(U^g,V^g,W^g)$ is the constant grid
velocity in a time interval. Due to the variation of the control
volume, the semi-discretized finite volume scheme can be expressed
as
\begin{align*}
\frac{\text{d}(|\Omega_{ijk}|Q_{ijk})}{\text{d}t}=-\mathcal{L}(Q_{ijk}),
\end{align*}
where the operator $\mathcal{L}$ is also given by Eq.\ref{finite}
and $|\Omega_{ijk}|$ varies in a time interval. In order to
update the flow variables in the moving framework, the numerical
fluxes at Gaussian quadrature points in Eq.\eqref{flux-gau} need to
be replaced by the following one with the mesh velocity
\begin{align*} F_{m_1,m_2}(t)=\int\psi f(\boldsymbol{x}_{m_1,m_2},t,\boldsymbol{u},\xi)
(\boldsymbol{u}-\boldsymbol{U}_{m_1,m_2}^g) \cdot
\boldsymbol{n}_{m_1,m_2}\text{d}\Xi,
\end{align*}
where $\boldsymbol{U}_{m_1,m_2}^g$ is the grid velocity at
quadrature point, which is given by the following interpolation
procedure
\begin{align*}
\boldsymbol{U}^g(\eta,\zeta)=\sum_{i=1}^4\boldsymbol{U}^g_i\phi_i(\eta,\zeta),
\end{align*}
where $\boldsymbol{U}^g_i$ is the velocity of four vertexes. With
the above procedure, the numerical fluxes and its temporal derivative at
$t_n$ are given as follows
\begin{align*}
\mathbb{F}_{p}^n&=\sum_{m_1,m_2=1}^2\omega_{m_1,m_2}
F_{m_1,m_2}^n\|\boldsymbol{X}_\eta\times\boldsymbol{X}_\zeta\|_{m_1,m_2}^n\Delta\eta\Delta\zeta,\\
\partial_t\mathbb{F}_{p}^n&=\sum_{m_1,m_2=1}^2\omega_{m_1,m_2}
\partial_tF_{m_1,m_2}^n\|\boldsymbol{X}_\eta\times\boldsymbol{X}_\zeta\|_{m_1,m_2}^n\Delta\eta\Delta\zeta,
\end{align*}
where
$\|\boldsymbol{X}_\eta\times\boldsymbol{X}_\zeta\|_{m_1,m_2}^n$
is the geometrical information at $t_n$. Similarly, the numerical
fluxes and temporal derivatives at intermediate state can be obtained
as well.
\end{rmk}

\section{Numerical tests}
In this section, numerical tests for both inviscid and viscous flows
will be presented to validate the current scheme. For the inviscid
flows, the collision time $\tau$ takes
\begin{align*}
\tau=\epsilon \Delta t+C\displaystyle|\frac{p_l-p_r}{p_l+p_r}|\Delta
t,
\end{align*}
where $\epsilon=0.01$ and $C=1$. For the viscous flows, we have
\begin{align*}
\tau=\frac{\mu}{p}+C \displaystyle|\frac{p_l-p_r}{p_l+p_r}|\Delta t,
\end{align*}
where $p_l$ and $p_r$ denote the pressure on the left and right
sides of the cell interface, $\mu$ is the dynamic viscous
coefficient, and $p$ is the pressure at the cell interface. In
smooth flow regions, it will reduce to $\tau=\mu/p$. Without special
statement, the specific heat ratio $\gamma=1.4$ and the CFL number
$CFL=0.35$ are used in the computation.

To improve the robustness, a simple limiting procedure is used. For
the reconstructed variables $P_m(\boldsymbol{x}_{G}), m=0,1,...,24$
from quadrature and linear polynomials, if any one value of the
densities $\rho_m(\boldsymbol{x}_{G})$ and pressures
$p_m(\boldsymbol{x}_{G}), m=0,1,...,24$ become negative, the
derivatives are set as zero and first-order reconstruction is
adopted. In order to eliminate the spurious oscillation and improve
the stability, the reconstruction can be performed for the
characteristic variables. The characteristic variables are defined
as $U=R^{-1}Q$, where $R$ is the right eigenmatrix of Jacobian
matrix $n_x(\partial F/\partial Q)_{G}+n_y(\partial G/\partial
Q)_{G}+n_z(\partial H/\partial Q)_{G}$ at Gaussian quadrature point.
With the reconstructed values, the conservative variables can be
obtained by the inverse projection.

\begin{table}[!h]
\begin{center}
\def\temptablewidth{0.75\textwidth}
{\rule{\temptablewidth}{1.0pt}}
\begin{tabular*}{\temptablewidth}{@{\extracolsep{\fill}}c|cc|cc}
mesh     & $L^1$ error  &    Order    &  $L^2$ error &  Order   \\
\hline
$16^3$   &    1.4612E-01  &    ~~      &  5.7820E-02     &   ~~    \\
$32^3$   &    2.0241E-02  &  2.8517    &  7.9355E-03     &  2.8651 \\
$64^3$   &    2.5712E-03  &  2.9768    &  1.0083E-03     &  2.9763 \\
$128^3$  &   3.2240E-04   &  2.9955    &  1.2633E-04     &  2.9966 \\
\end{tabular*}
{\rule{\temptablewidth}{1.0pt}}
\end{center}
\vspace{-2.5mm}\caption{\label{tab-3d-1} Accuracy test: the
advection of density perturbation with uniform meshes.}
\begin{center}
\def\temptablewidth{0.75\textwidth}
{\rule{\temptablewidth}{1.0pt}}
\begin{tabular*}{\temptablewidth}{@{\extracolsep{\fill}}c|cc|cc}
mesh     & $L^1$ error  &    Order      &  $L^2$ error &  Order   \\
\hline
$16^3$   &  1.9151E-01    &   ~~        &  7.5334E-02    &    ~~   \\
$32^3$   &  2.8252E-02    &  2.7610     &  1.1099E-02    &   2.7628\\
$64^3$   &  3.6640E-03    &  2.9468     &  1.4352E-03    &   2.9511\\
$128^3$  &  4.6186E-04    &  2.9879     &  1.8063E-04   &   2.9901\\
\end{tabular*}
{\rule{\temptablewidth}{1.0pt}}
\end{center}
\vspace{-2.5mm}\caption{\label{tab-3d-2} Accuracy test: the
advection of density perturbation with non-coplanar meshes.}
\end{table}

\subsection{Accuracy test}
The advection of density perturbation for three-dimensional flows is
presented to test the order of accuracy. For this case, the physical
domain is $[0,2]\times[0,2]\times[0,2]$ and the initial condition is
set as follows
\begin{align*}
\rho_0&(x, y, z)=1+0.2\sin(\pi(x+y+z)),~p_0(x,y,z)=1,\\
&U_0(x,y,z)=1,~V_0(x,y,z)=1,~W_0(x,y,z)=1.
\end{align*}
The periodic boundary conditions are applied at boundaries, and the
exact solution is
\begin{align*}
\rho(x,y&,z,t)=1+0.2\sin(\pi(x+y+z-t)),~p(x,y,z,t)=1,\\
&U(x,y,z,t)=1,~V(x,y,z,t)=1,~W(x,y,z,t)=1.
\end{align*}
The uniform mesh with $\Delta x=\Delta y=\Delta z=2/N$ are tested.
The $L^1$ and $L^2$ errors and orders of accuracy at $t=2$ are
presented in Tab.\ref{tab-3d-1}, where the expected order of
accuracy is achieved. To validate the order of accuracy with
non-coplanar meshes, the following mesh is considered
\begin{align*}
\begin{cases}
\displaystyle x_i=\xi_i+0.1\sin (\pi \xi_i)\sin (\pi \eta_j)\sin (\pi \zeta_k),\\
\displaystyle y_j=\eta_j+0.1\sin (\pi \xi_i)\sin (\pi \eta_j)\sin (\pi \zeta_k),\\
\displaystyle z_k=\zeta_k+0.1\sin (\pi \xi_i)\sin (\pi \eta_j)\sin
(\pi \zeta_k),
\end{cases}
\end{align*}
where $(\xi,\eta,\zeta)\in [0,2]\times[0,2]\times[0,2]$, and
$(\xi_i,\eta_j,\zeta_k)$ are given uniformly with $\Delta \xi=\Delta
\eta=\Delta \zeta=2/N$. For most cells given above, it can be easily
verified that the vertexes are non-coplanar. The $L^1$ and $L^2$
errors and orders of accuracy at $t=2$ are presented in
Tab.\ref{tab-3d-2}, and the expected order of accuracy is achieved
by the current scheme as well. For the mesh obtained by smooth
coordinate transformation, numerical scheme can be constructed with
structured WENO reconstruction, and more details can be found in
\cite{GKS-high-2}.

\subsection{Moving-mesh tests}
To validate the order of accuracy with moving-meshes, the following
mesh is considered
\begin{align*}
\begin{cases}
\displaystyle x_i=\xi_i+0.1\sin (\pi \xi_i)\sin (\pi \eta_j)\sin (\pi \zeta_k)\sin\pi t,\\
\displaystyle y_j=\eta_j+0.1\sin (\pi \xi_i)\sin (\pi \eta_j)\sin (\pi \zeta_k)\sin\pi t,\\
\displaystyle z_k=\zeta_k+0.1\sin (\pi \xi_i)\sin (\pi \eta_j)\sin
(\pi \zeta_k)\sin\pi t,
\end{cases}
\end{align*}
where $(\xi,\eta,\zeta)\in [0,2]\times[0,2]\times[0,2]$, and
$(\xi_i,\eta_j,\zeta_k)$ are given uniformly with $\Delta \xi=\Delta
\eta=\Delta \zeta=2/N$. The periodic boundary condition is imposed
for the mesh. The $L^1$ and $L^2$ errors and orders of accuracy
after one period, i.e. $t=2$ are presented in Tab.\ref{tab-3d-3}.
The order of accuracy is well kept during the moving-mesh procedure.
The geometric conservation law \cite{GCL} is also tested, which is
mainly about the maintenance of uniform flow passing through the
moving-mesh. The initial condition is given as follows
\begin{align*}
\rho_0&(x, y,
z)=1,~p_0(x,y,z)=1,~U_0(x,y,z)=1,~V_0(x,y,z)=1,~W_0(x,y,z)=1.
\end{align*}
The above moment of computational mesh is used, and the periodic
boundary conditions are adopted. The $L^1$ and $L^2$ errors at
$t=0.5$ are given in Tab.\ref{tab-3d-4}. The results show that the
errors reduce to the machine zero, which implies the satisfaction of
geometric conservation law.

\begin{table}[!h]
\begin{center}
\def\temptablewidth{0.75\textwidth}
{\rule{\temptablewidth}{1.0pt}}
\begin{tabular*}{\temptablewidth}{@{\extracolsep{\fill}}c|cc|cc}
mesh     & $L^1$ error  & Order &    $L^2$ error &  Order     \\
\hline
$16^3$   &1.6475E-01    & ~~    &    6.4872E-02  &    ~~       \\
$32^3$   &2.3757E-02    &2.7938 &    9.4385E-03  &    2.7809   \\
$64^3$   &3.0657E-03    &2.9540 &    1.2195E-03  &    2.9521   \\
$128^3$  &3.8585E-04    &2.9901 &    1.5341E-04  &    2.9908   \\
\end{tabular*}
{\rule{\temptablewidth}{1.0pt}}
\end{center}
\vspace{-2.5mm}\caption{\label{tab-3d-3} Moving-mesh test: 3D
advection of density perturbation with moving-meshes.}
\begin{center}
\def\temptablewidth{0.5\textwidth}
{\rule{\temptablewidth}{1.0pt}}
\begin{tabular*}{\temptablewidth}{@{\extracolsep{\fill}}c|cc}
3D   mesh   & $L^1$ error  &  $L^2$ error   \\
\hline
$16^3$   &   1.1925E-14 & 5.4002E-15\\
$32^3$   &   3.1108E-14 & 1.4144E-14\\
$64^3$   &   7.7908E-14 & 3.6674E-14\\
$128^3$  &   1.8682E-13 & 8.9173E-14
\end{tabular*}
{\rule{\temptablewidth}{1.0pt}}
\end{center}
\vspace{-2.5mm}\caption{\label{tab-3d-4} Moving-mesh test: the
geometric conservation law  with moving-meshes.}
\end{table}

\begin{figure}[!h]
\centering
\includegraphics[width=0.45\textwidth]{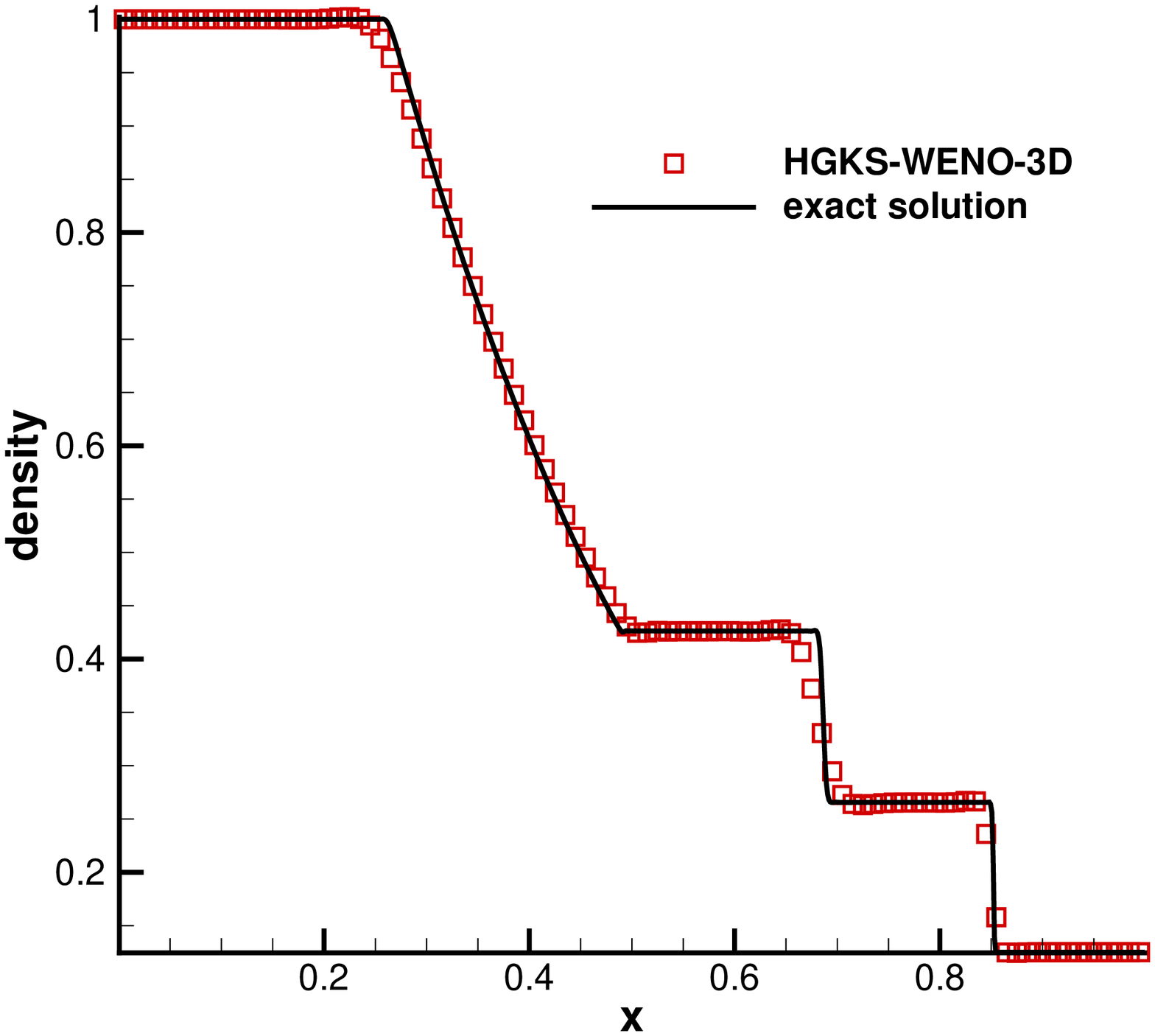}\includegraphics[width=0.45\textwidth]{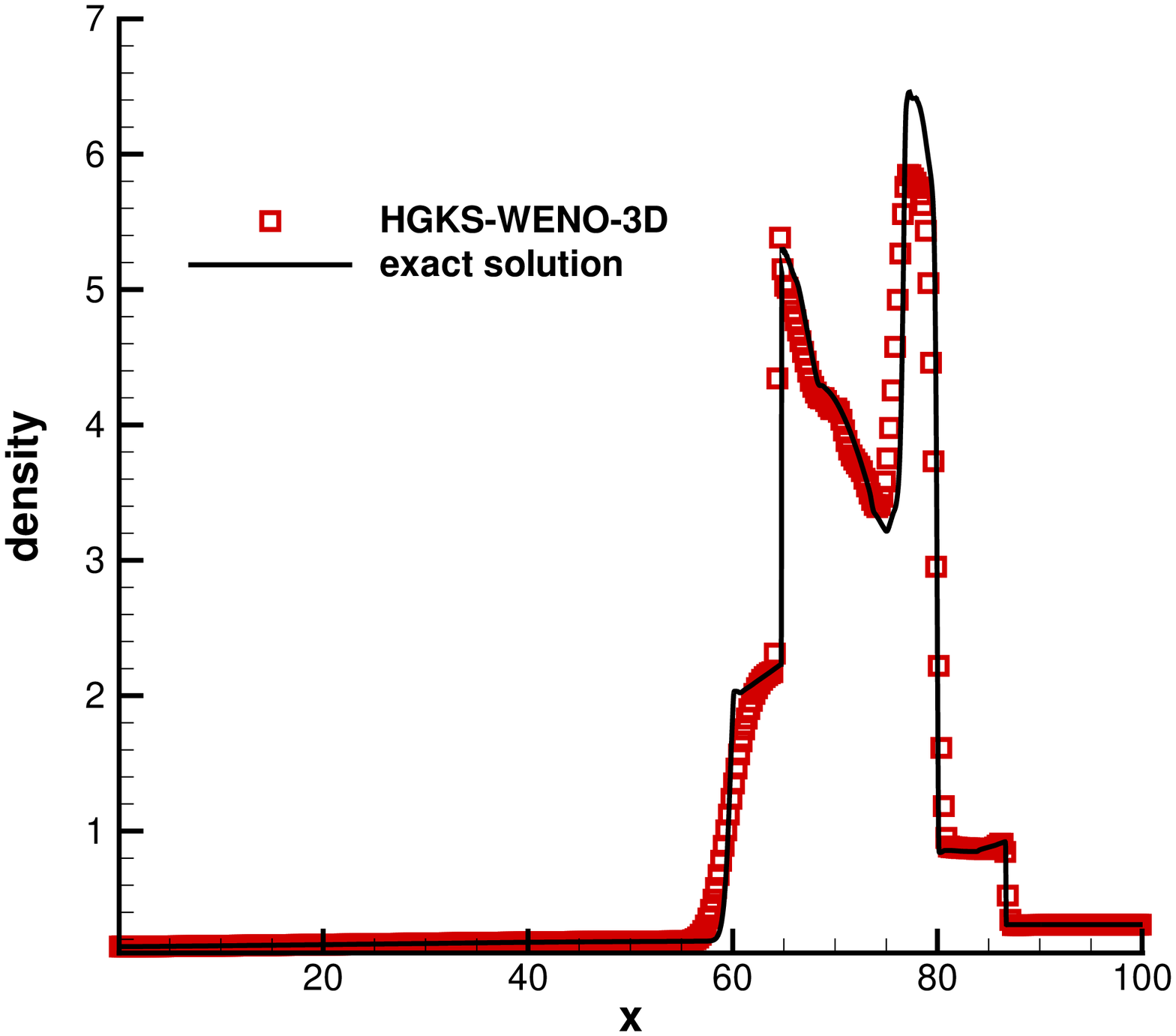}\\
\includegraphics[width=0.45\textwidth]{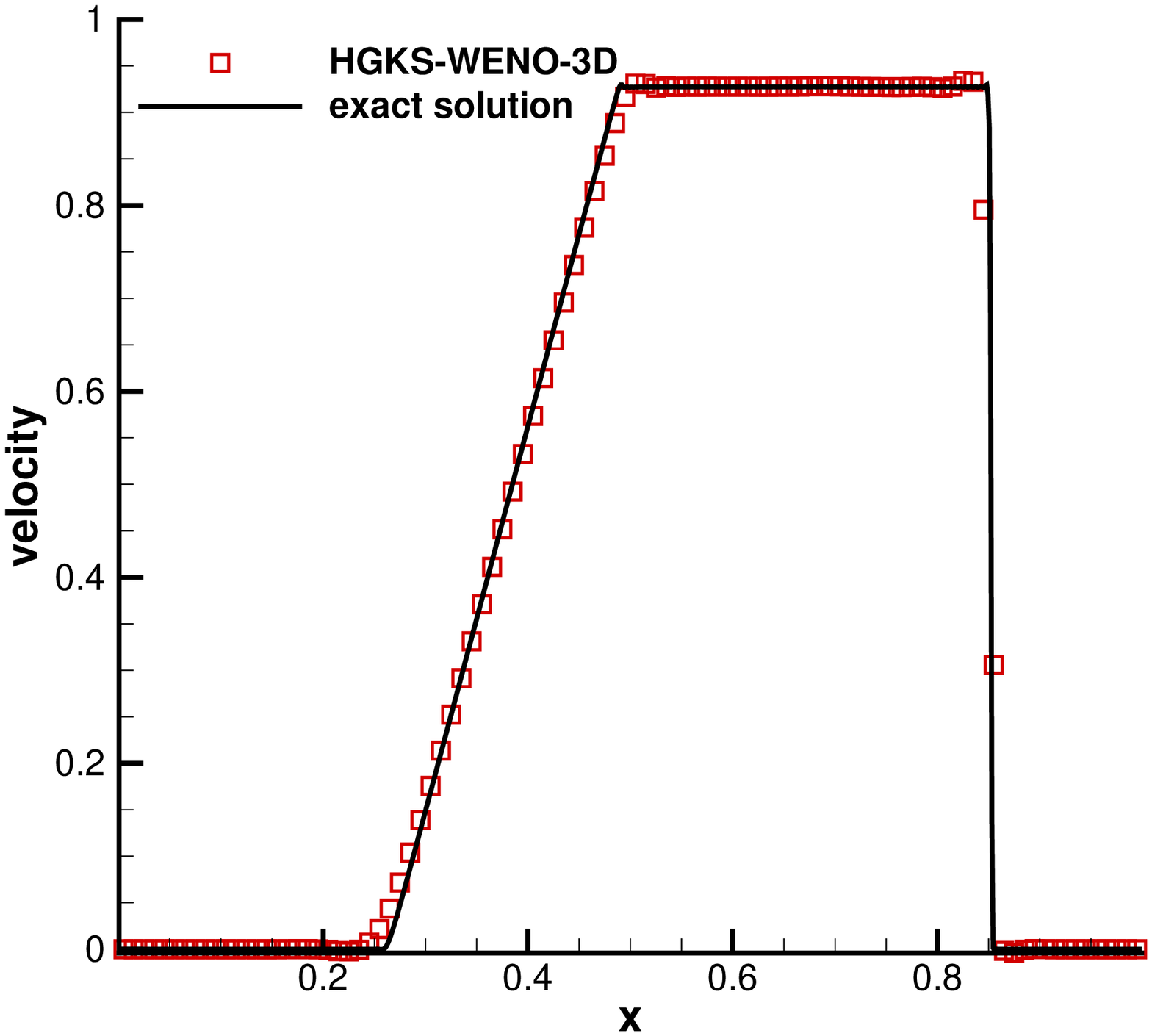}\includegraphics[width=0.45\textwidth]{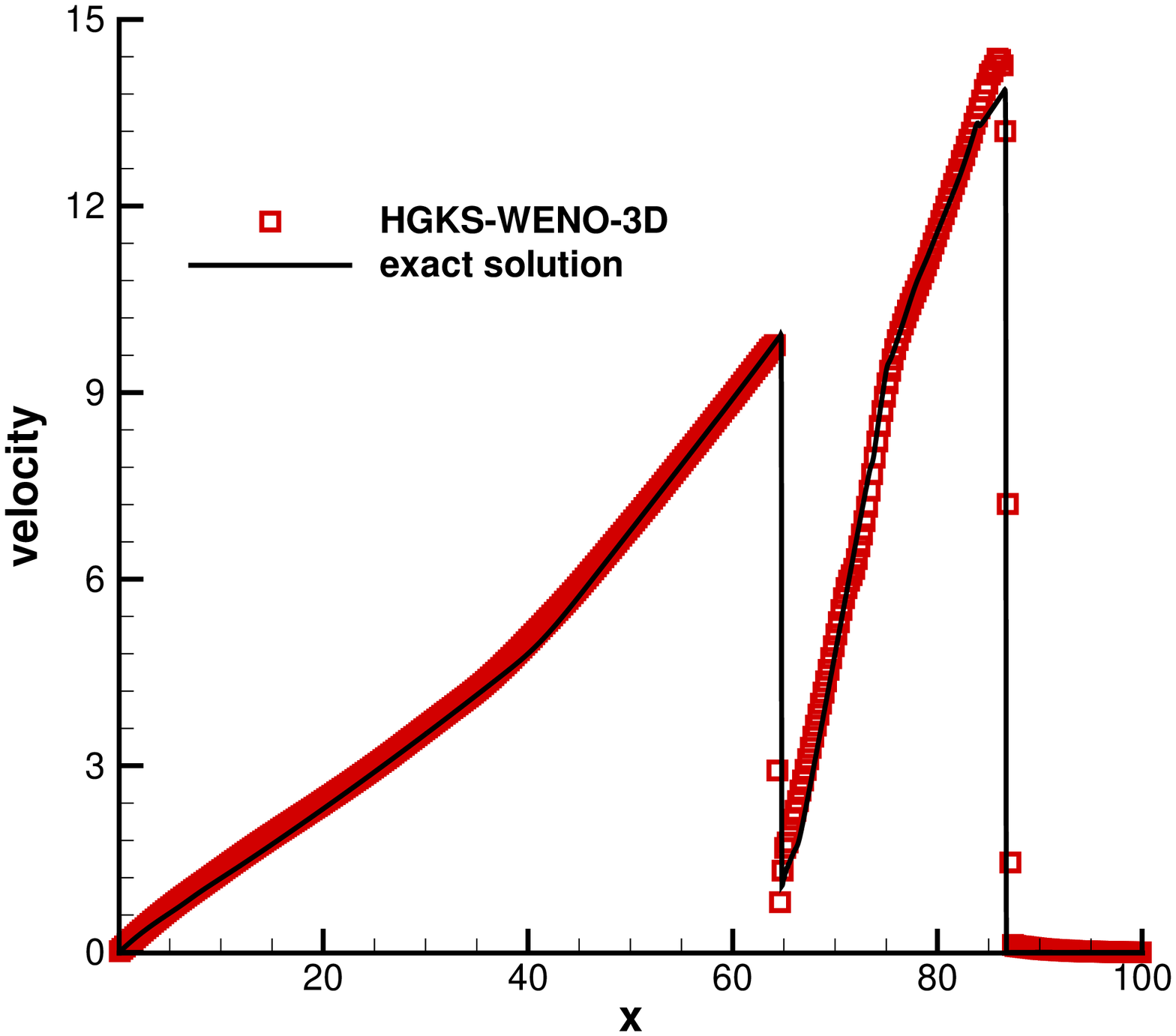}\\
\includegraphics[width=0.45\textwidth]{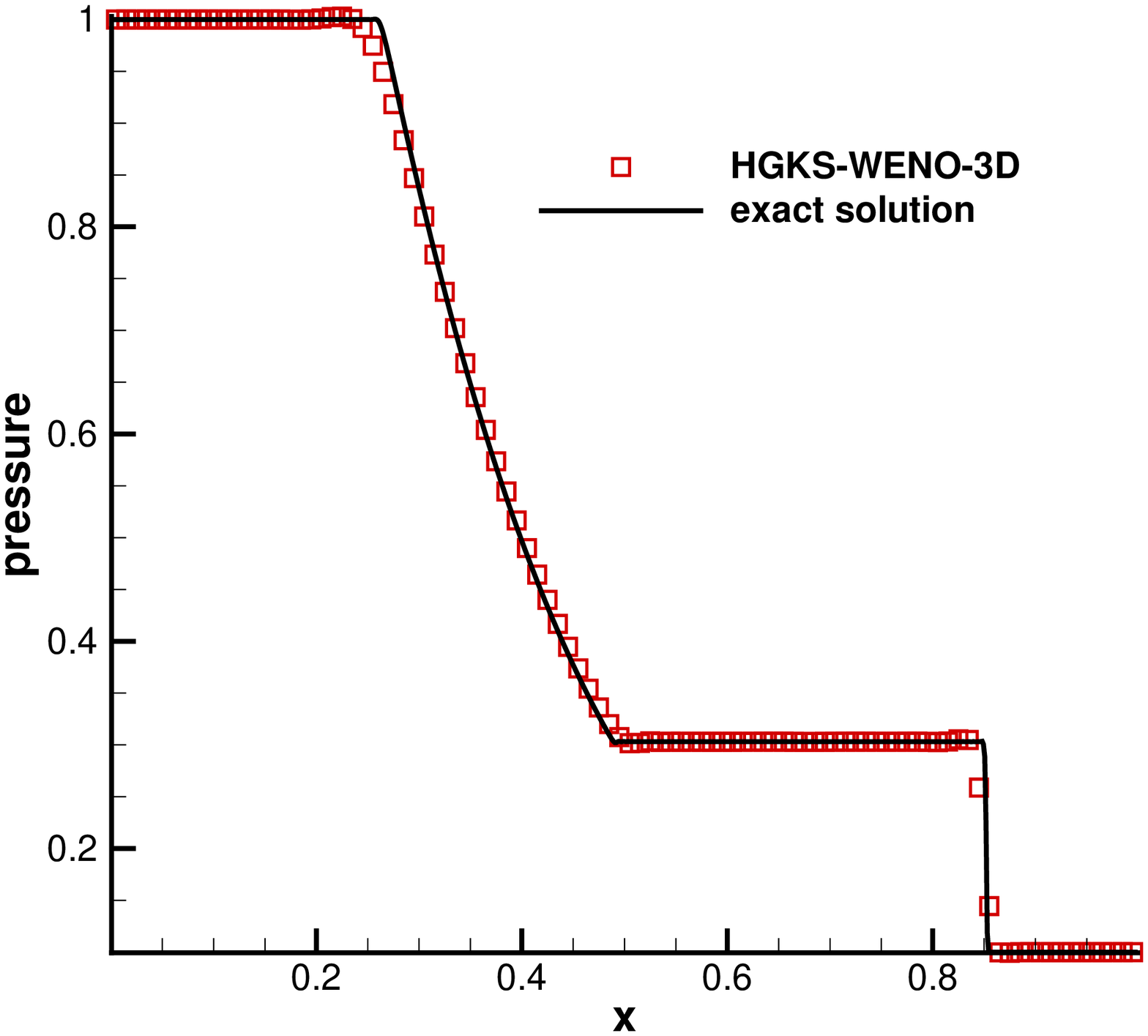}\includegraphics[width=0.45\textwidth]{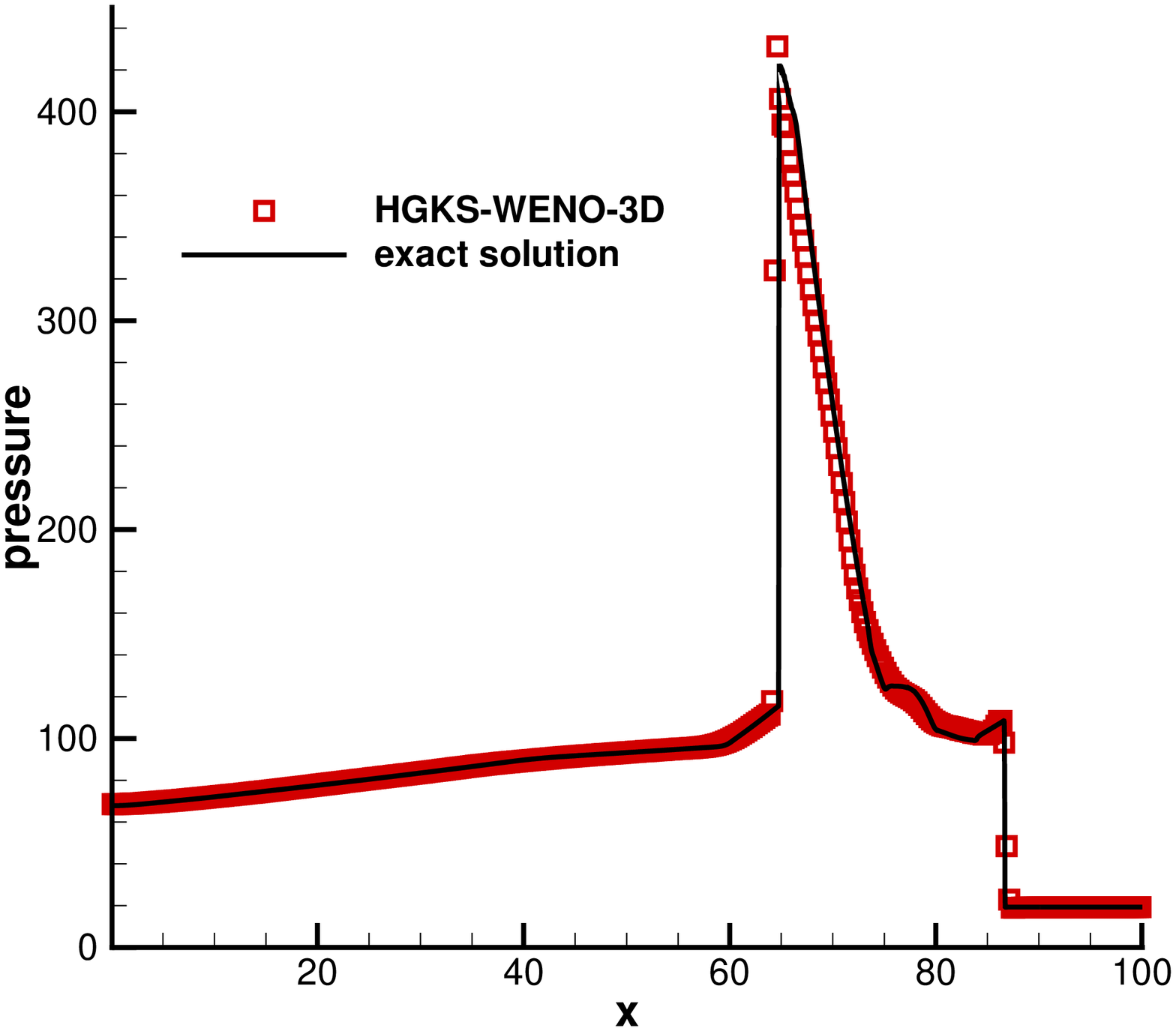}
\caption{\label{riemann-1} Riemann problem: the density, velocity
and pressure distributions at $t=0.2$ for Sod problem (left), and the density, velocity and pressure distributions at
$t=3.8$ for blast wave problem
(right).}
\end{figure}

\subsection{Riemann problem}
In this case, one-dimensional Riemann problems are considered. The
first one is the Sod problem, and the initial conditions are given
by
\begin{equation*}
(\rho, U, V, W, p)=
\begin{cases}
(1, 0, 0, 0, 1),  \ \ \ \ &  0\leq x<0.5,\\
(0.125, 0, 0, 0, 0.1),    & 0.5<x\leq1.
\end{cases}
\end{equation*}
The computational domain is $(x,y,z)\in[0, 1]\times[0, 0.1]\times[0,
0.1]$, and the uniform mesh with $\Delta x=\Delta y=\Delta z=1/100$
is used. The non-reflected boundary condition is used in all
directions. The second one is the Woodward-Colella blast wave
problem, which is used to test the robustness of WENO
reconstruction. The initial conditions are given as follows,
\begin{equation*}
(\rho, U, V, W, p) =\left\{\begin{array}{ll}
(1, 0, 0, 0, 1000), \ \ \ \ & 0\leq x<10,\\
(1, 0, 0, 0, 0.01), & 10\leq x<90,\\
(1, 0, 0, 0, 100), &  90\leq x\leq 100.
\end{array} \right.
\end{equation*}
The computational domain is $(x,y,z)\in[0, 100]\times[0,
10]\times[0,10]$, and the uniform mesh with $\Delta x=\Delta
y=\Delta z=1/4$ is used.  The reflected boundary condition is used
in $x$ direction, and non-reflected boundary condition is used in
$y$ and $z$ directions. The density, velocity, and pressure
distributions for the current scheme and the exact solutions are
presented in Fig.\ref{riemann-1} for the Sod problem at $t=0.2$ and
for the blast wave problem at $t=3.8$ with  $y=z=0$. The numerical
results agree well with the exact solutions.

\begin{figure}[!h]
\centering
\includegraphics[width=0.45\textwidth]{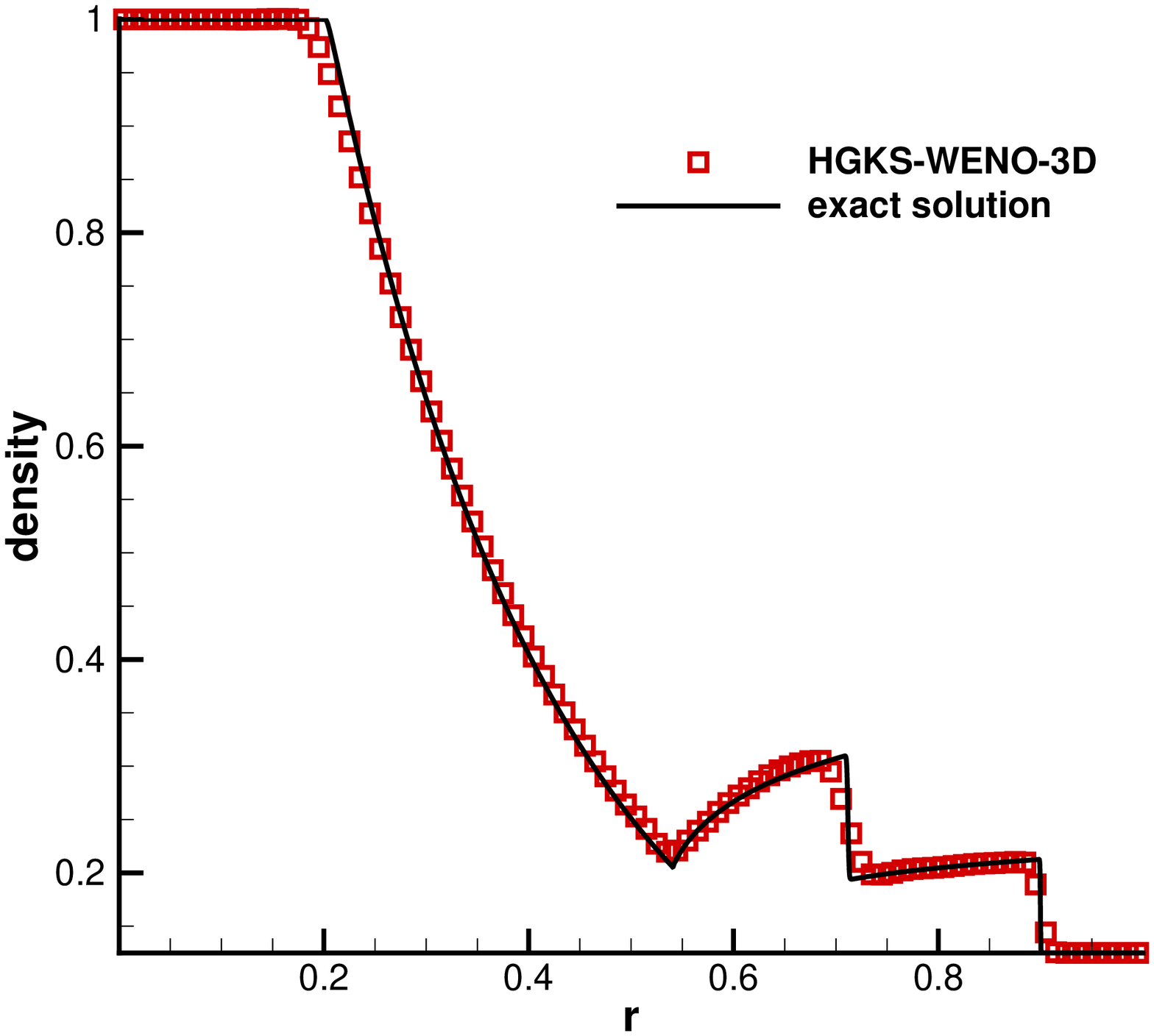}
\includegraphics[width=0.45\textwidth]{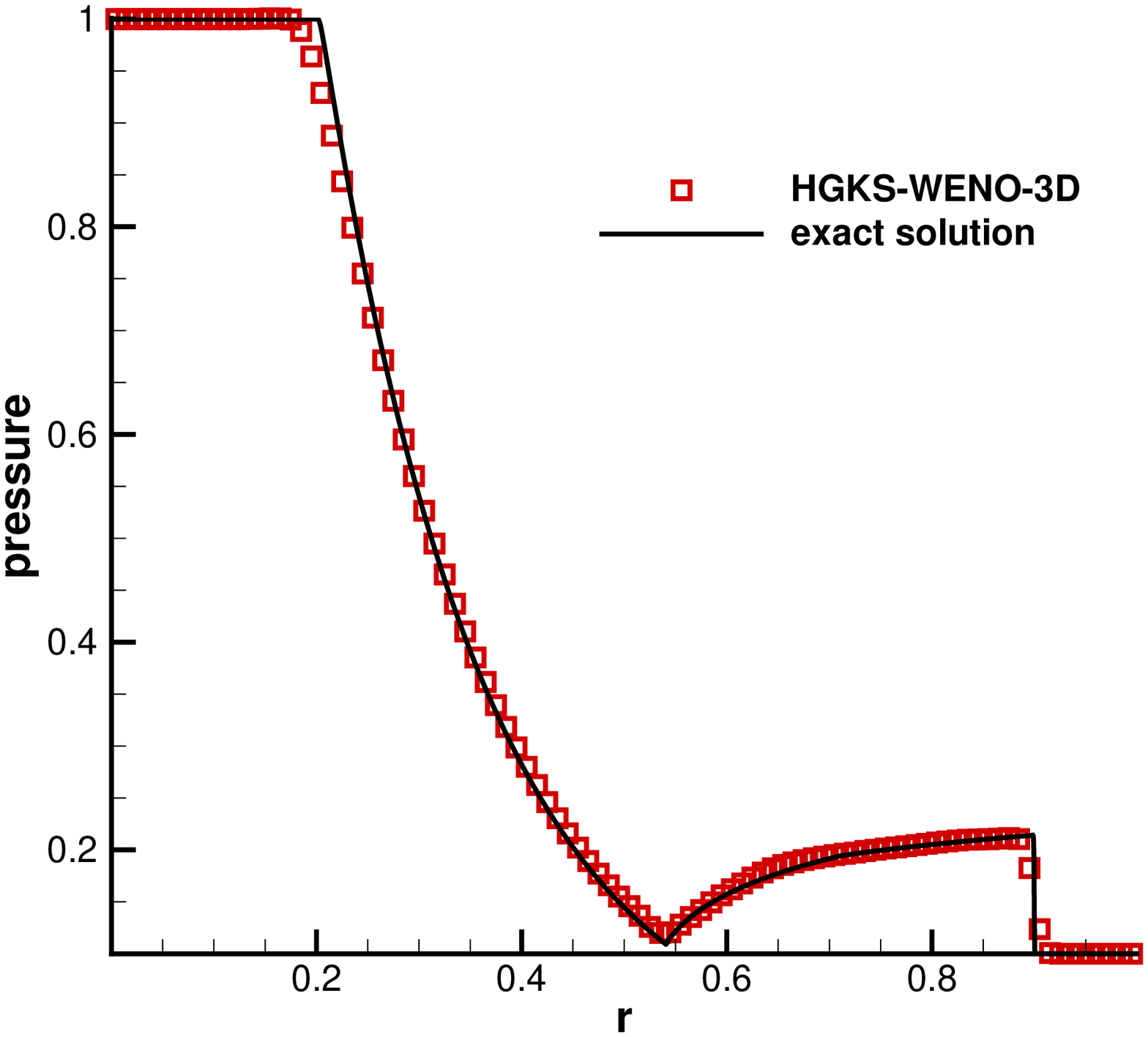}
\caption{\label{riemann-2} Riemann problem: the density and pressure
distributions at $t=0.25$ for spherically symmetric Sod problem .}
\end{figure}

As an extension of the Sod problem, the spherically symmetric Sod
problem is considered, and the initial conditions are given by
\begin{equation*}
(\rho, U, V, W, p)=
\begin{cases}
(1, 0, 0, 0, 1),  \ \ \ \ &  0\leq\sqrt{x^2+y^2+z^2}<0.5,\\
(0.125, 0, 0, 0, 0.1),    & 0.5<\sqrt{x^2+y^2+z^2}\leq1.
\end{cases}
\end{equation*}
The computational domain is $(x,y,z)\in[0, 1]\times[0, 1]\times[0,
1]$, and the uniform mesh with $\Delta x=\Delta y=\Delta z=1/100$ is
used. The symmetric boundary condition is imposed on the plane with
$x=0,y=0,z=0$, and the non-reflection boundary condition is imposed
on the plane with $x=1,y=1,z=1$. The exact solution of spherically
symmetric problem can be given by the following one-dimensional
system with geometric source terms
\begin{align*}
\frac{\partial Q}{\partial t}+\frac{\partial F(Q)}{\partial r}=S(Q),
\end{align*}
where
\begin{align*}
Q=\begin{pmatrix}
   \rho \\
   \rho U \\
   \rho E \\
\end{pmatrix},
F(Q)=\begin{pmatrix}
   \rho U\\
   \rho U^2+p \\
   U(\rho E+p) \\
 \end{pmatrix},
S(Q)=-\frac{d-1}{r}\begin{pmatrix}
   \rho U\\
   \rho U^2  \\
   U(\rho E+p) \\
 \end{pmatrix}.
\end{align*}
The radial direction is denoted by $r$, $U$ is the radial velocity,
$d$ is the number of space dimensions. The density and pressure
profiles along $y=z=0$ at $t=0.25$ are given in Fig.\ref{riemann-2}.
The current scheme also well resolves the wave profiles.

\subsection{Sedov blast wave problem}
\begin{figure}[!h]
\centering
\includegraphics[width=0.475\textwidth]{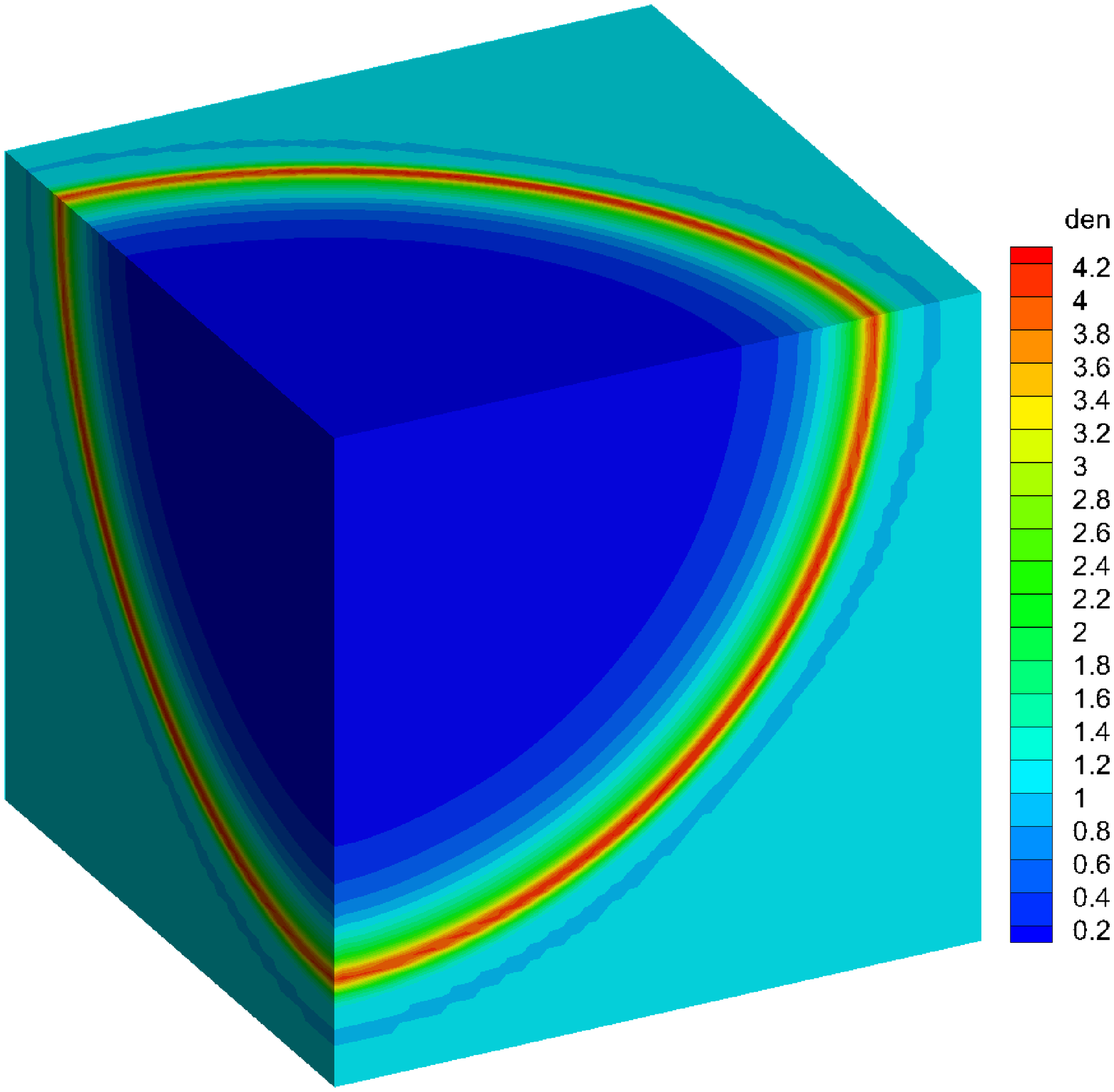}
\includegraphics[width=0.475\textwidth]{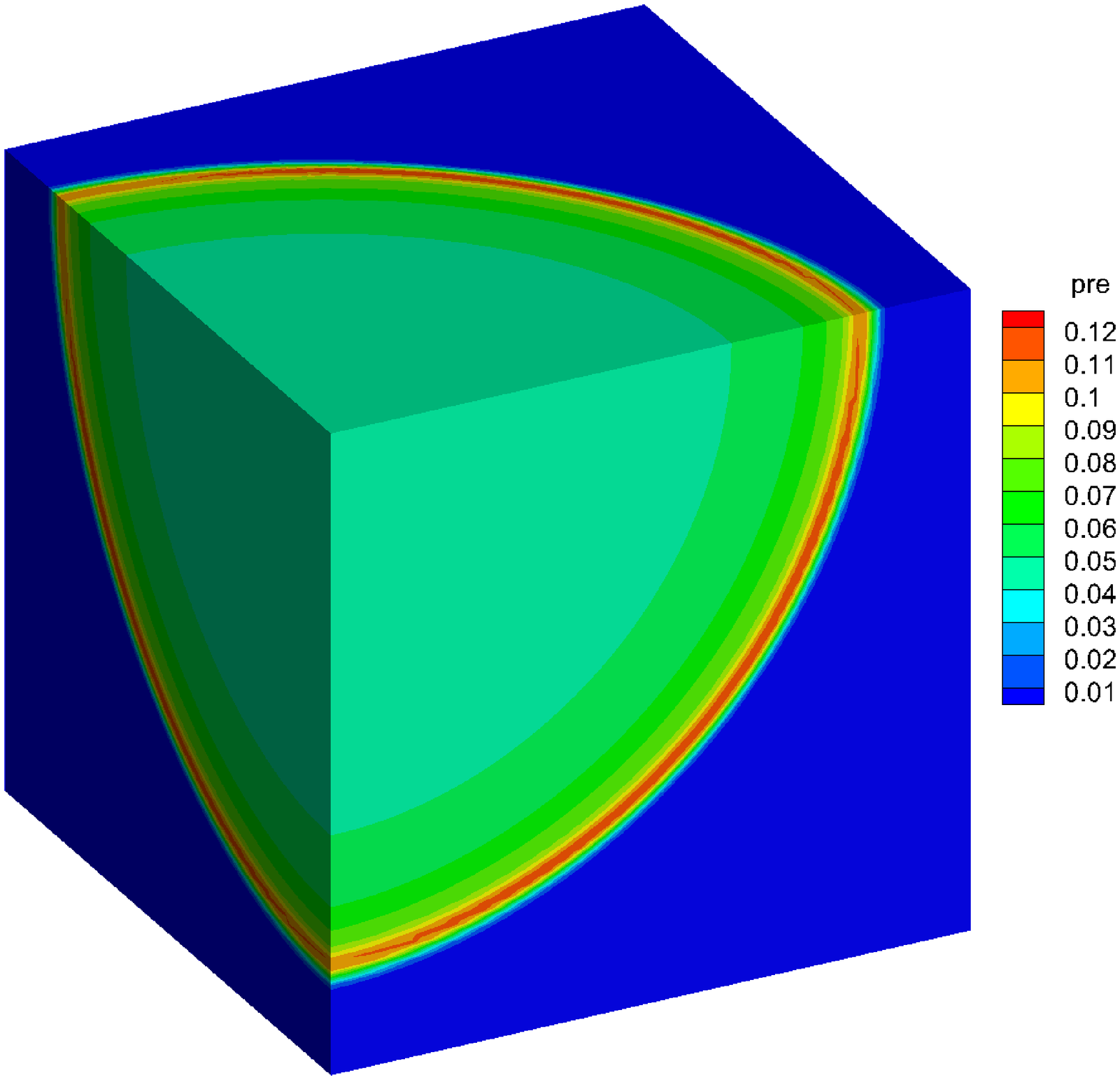}
\caption{\label{Sedov-a} Sedov problem:  the three-dimensional
density and pressure distributions at $t=1$ with $80\times 80\times
80$ cells.} \centering
\includegraphics[width=0.475\textwidth]{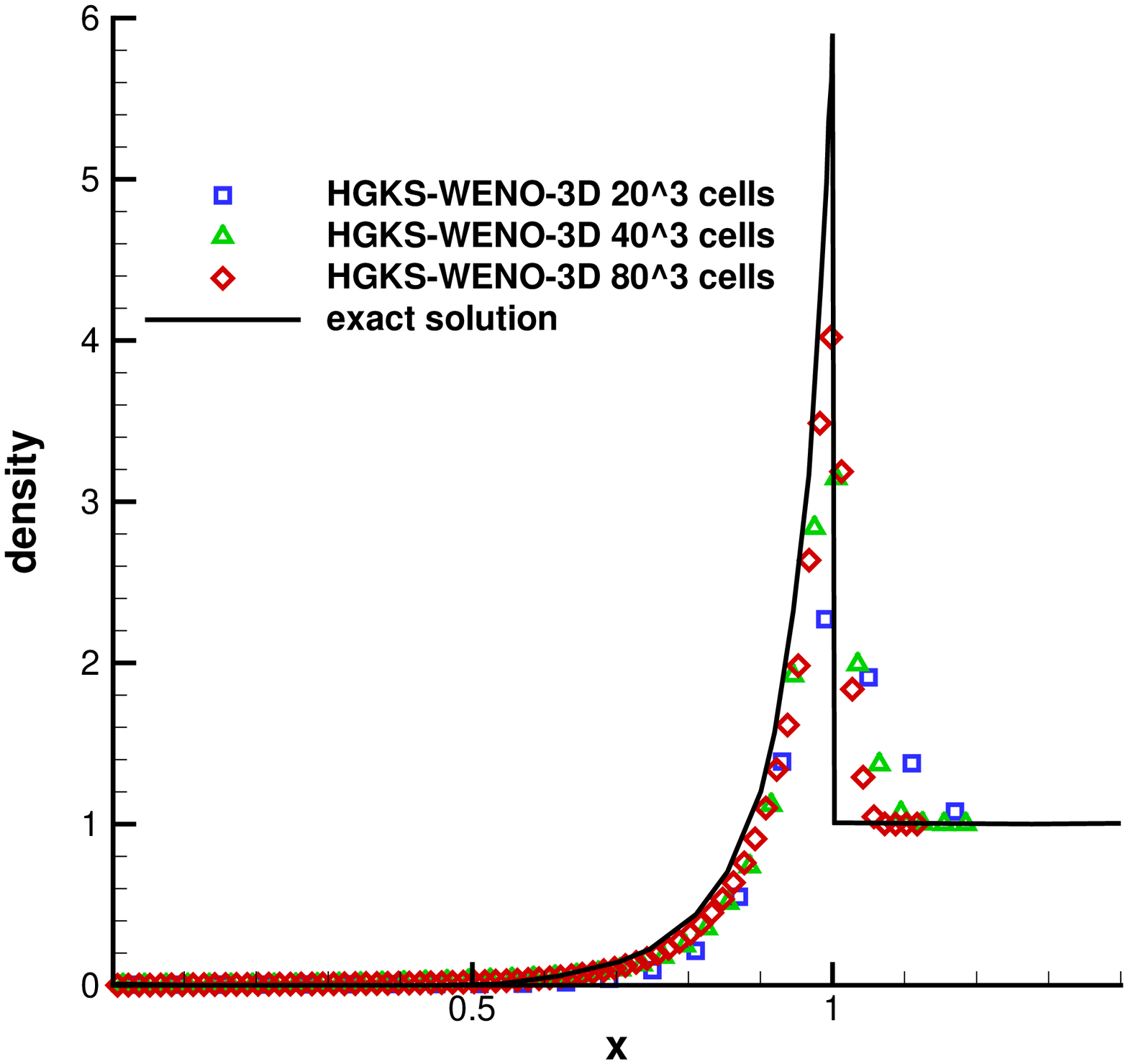}
\includegraphics[width=0.475\textwidth]{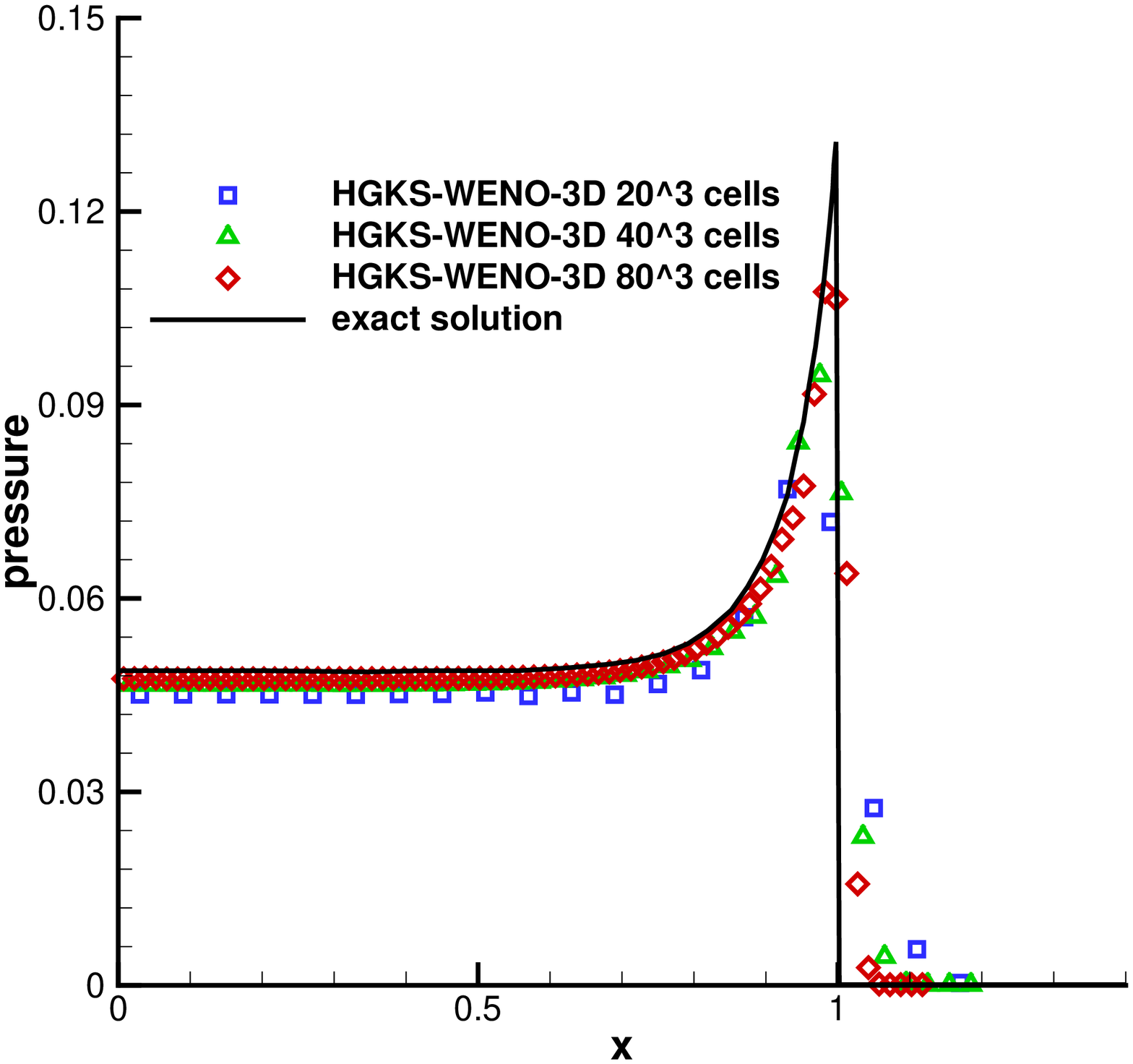}
\caption{\label{Sedov-b} Sedov problem: the density and pressure profiles
long $y=z=0$ at $t=1$.}
\end{figure}

This is a three-dimensional explosion problem to model blast wave
from an energy deposited singular point.  The initial density has a
uniform unit distribution, and the pressure is $10^{-6}$ everywhere,
except in the cell containing the origin. For this cell containing
the origin, the pressure is defined as
$p=(\gamma-1)\varepsilon_0/V$, where $\varepsilon_0=0.106384$ is the
total amount of released energy and $V$ is the cell volume. The
computation domain is $[0,1.2]\times[0,1.2]\times[0,1.2]$, and
uniform meshes are used. Due to the singularity at the origin, a
small CFL number $CFL=0.01$ is used. After 10 steps, a normal CFL
number is used. The symmetric boundary condition is imposed on the
plane with $x=0,y=0,z=0$, and the non-reflection boundary condition
is imposed on the plane with $x=1.2,y=1.2,z=1.2$. The solution
consists of a diverging infinite strength shock wave whose front is
located at radius $r=1$ at $t=1$ \cite{Case-Sedov}. The
three-dimensional density and pressure distributions with $80\times
80\times 80$ cells at $t=1$ are presented in Fig.\ref{Sedov-a}. The
density and pressure profiles long $y=z=0$ at $t=1$ with $20\times
20\times 20, 40\times 40\times 40$ and $80\times 80\times 80$ cells
are given in Fig.\ref{Sedov-b}. With the mesh refinement, the
numerical solutions agree well with the exact solutions.

\begin{figure}[!h]
\centering
\includegraphics[width=0.95\textwidth]{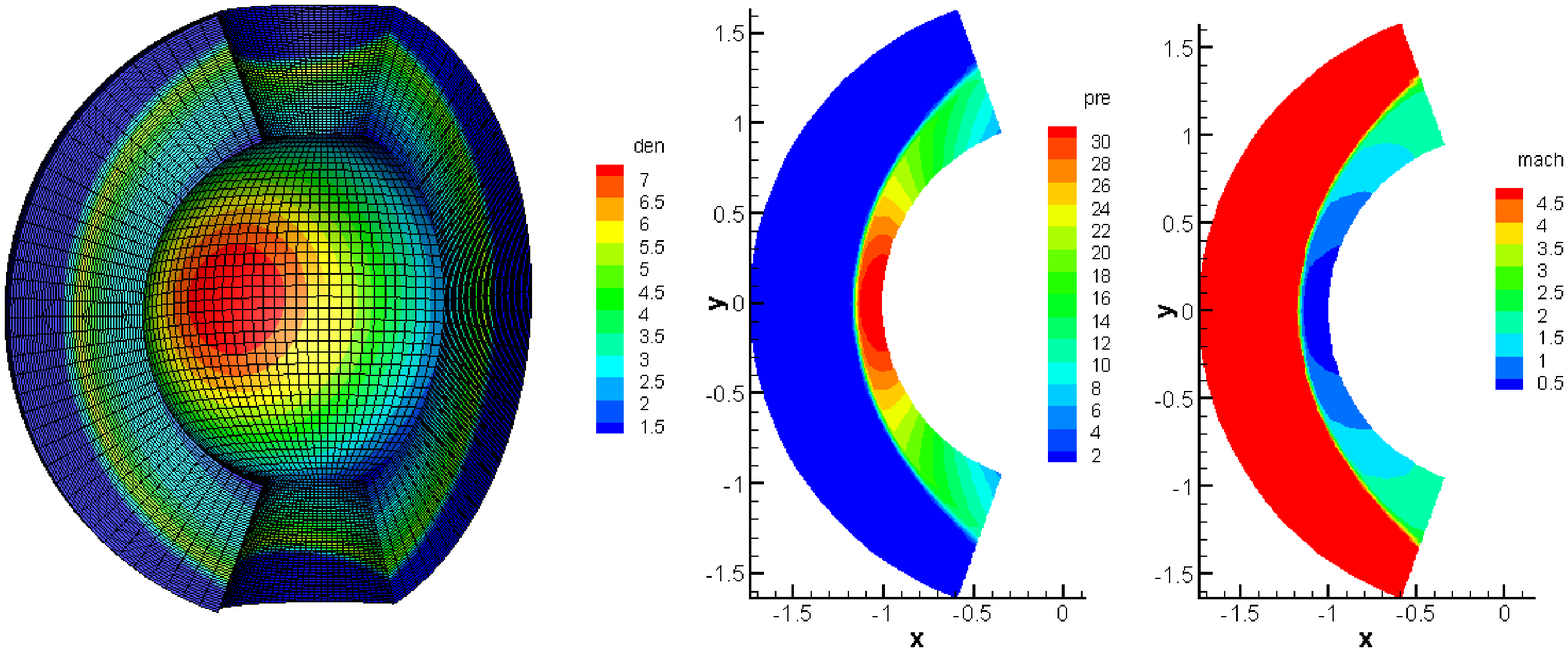}{a}
\includegraphics[width=0.95\textwidth]{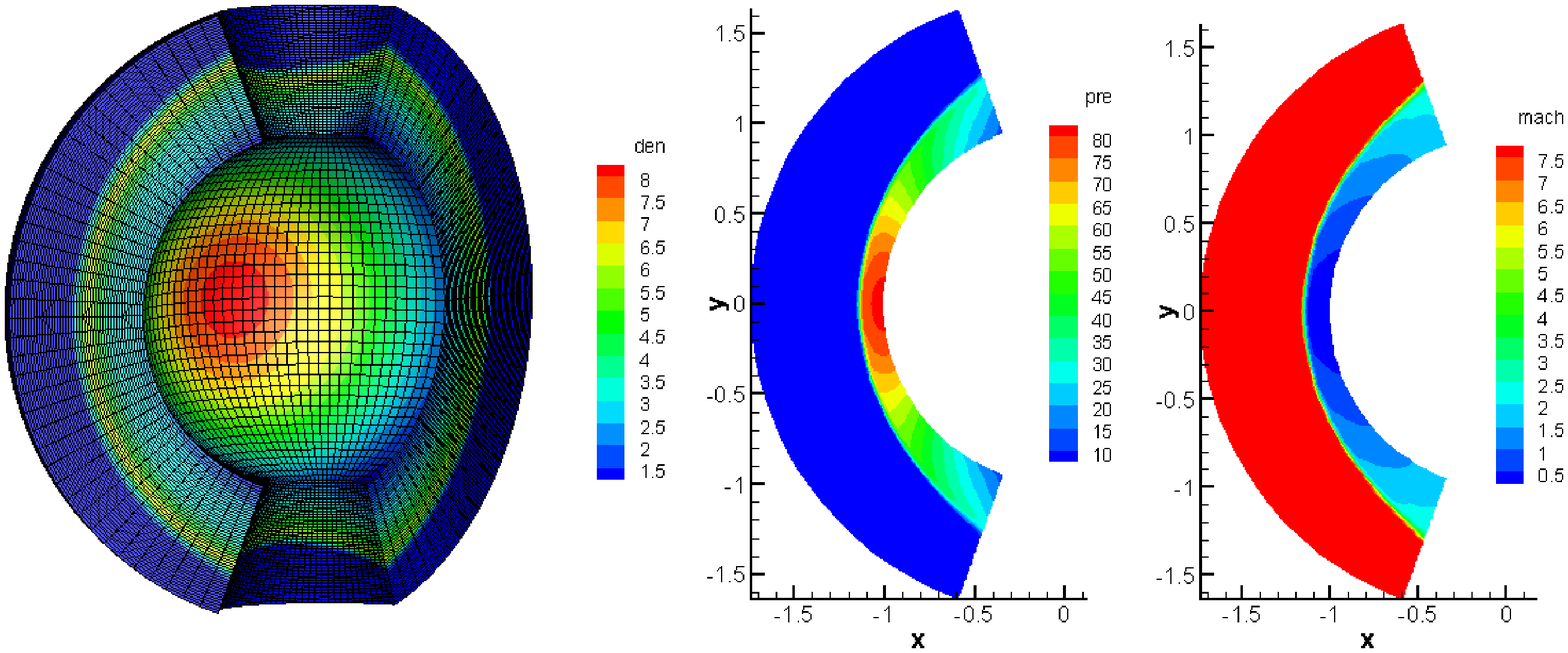}{b}
\caption{\label{sphere} Flow impinging on sphere: the mesh and
density distribution, pressure and Mach number distribution at
$\phi=0$ for $Ma=5$ (a) and $8$ (b).}
\end{figure}

\subsection{Flow impinging on sphere}
In this case, the inviscid hypersonic flows impinging on a unit
sphere are tested to validate robustness of the current scheme with
high Mach numbers. In the computation, a $50\times40\times40$ mesh
is used, which represents the domain
$[-1.75,-1]\times[-0.4\pi,0.4\pi]\times[0.1\pi,0.9\pi]$ in the
spherical coordinate $(r,\phi,\theta)$.  The mesh and density
distributions in the whole domain, pressure and Mach number
distributions at the plane with $\phi=0$ for the Mach number $Ma=5$
and $8$ are shown in Fig.\ref{sphere}, where the shock is well
captured by the current scheme and the carbuncle phenomenon does not
appear.

\begin{figure}[!h]
\centering
\includegraphics[width=0.475\textwidth]{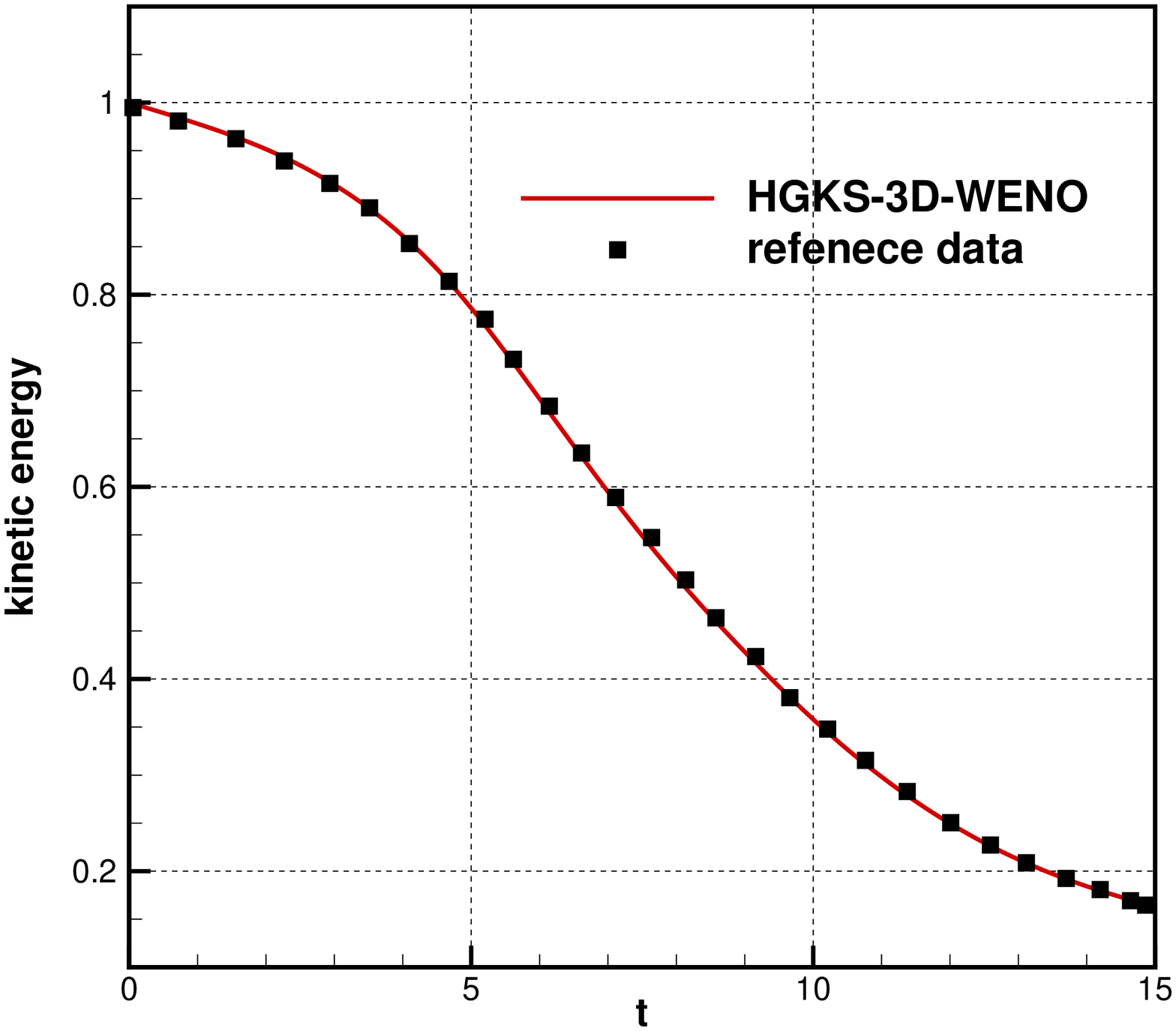}
\includegraphics[width=0.475\textwidth]{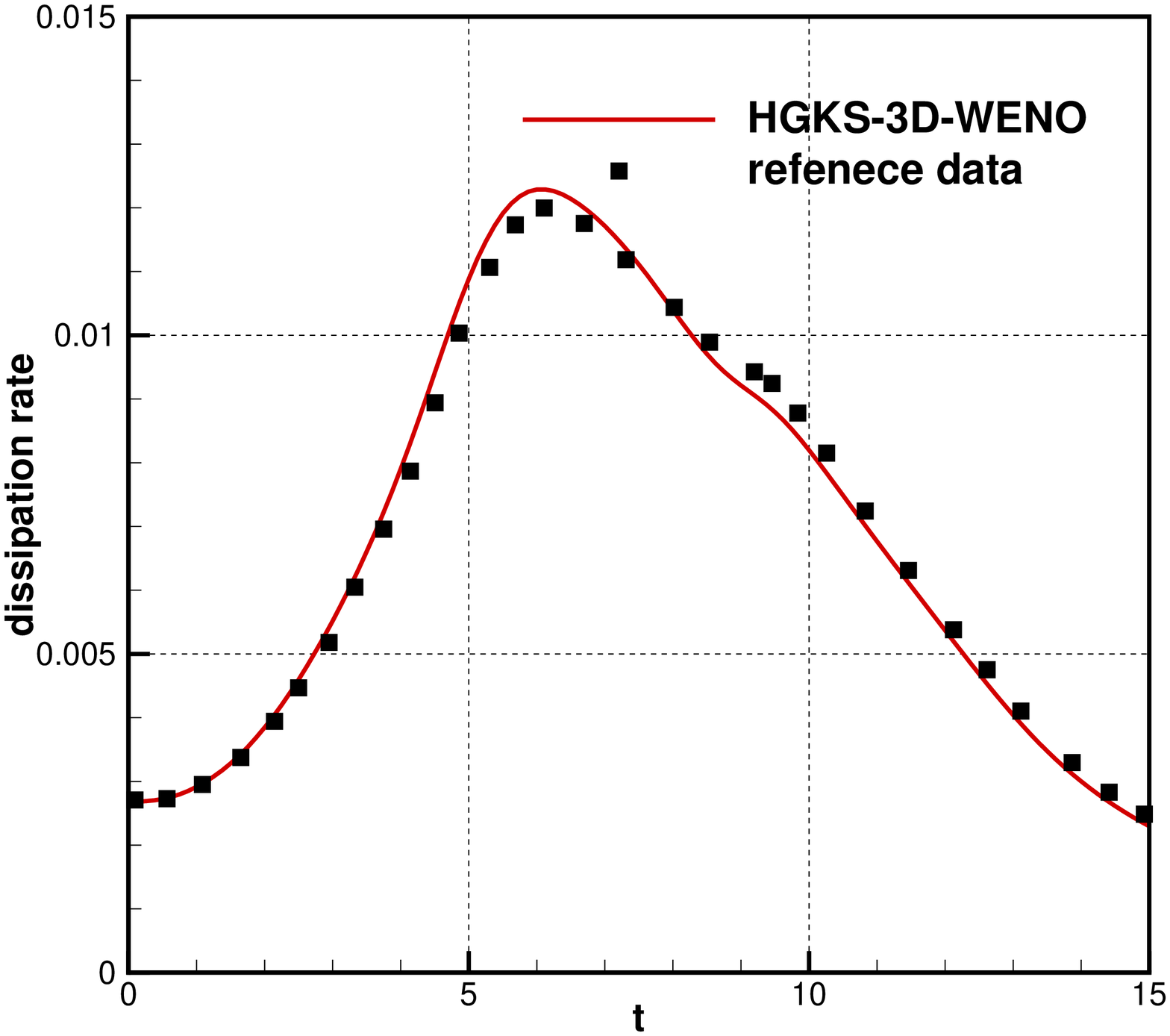}
\caption{\label{tg-vortex-1} Taylor-Green Vortex: the time history
of kinetic energy $E_k$ and dissipation rate
$-\text{d}E_k/\text{d}t$.}
\centering
\includegraphics[width=0.45\textwidth]{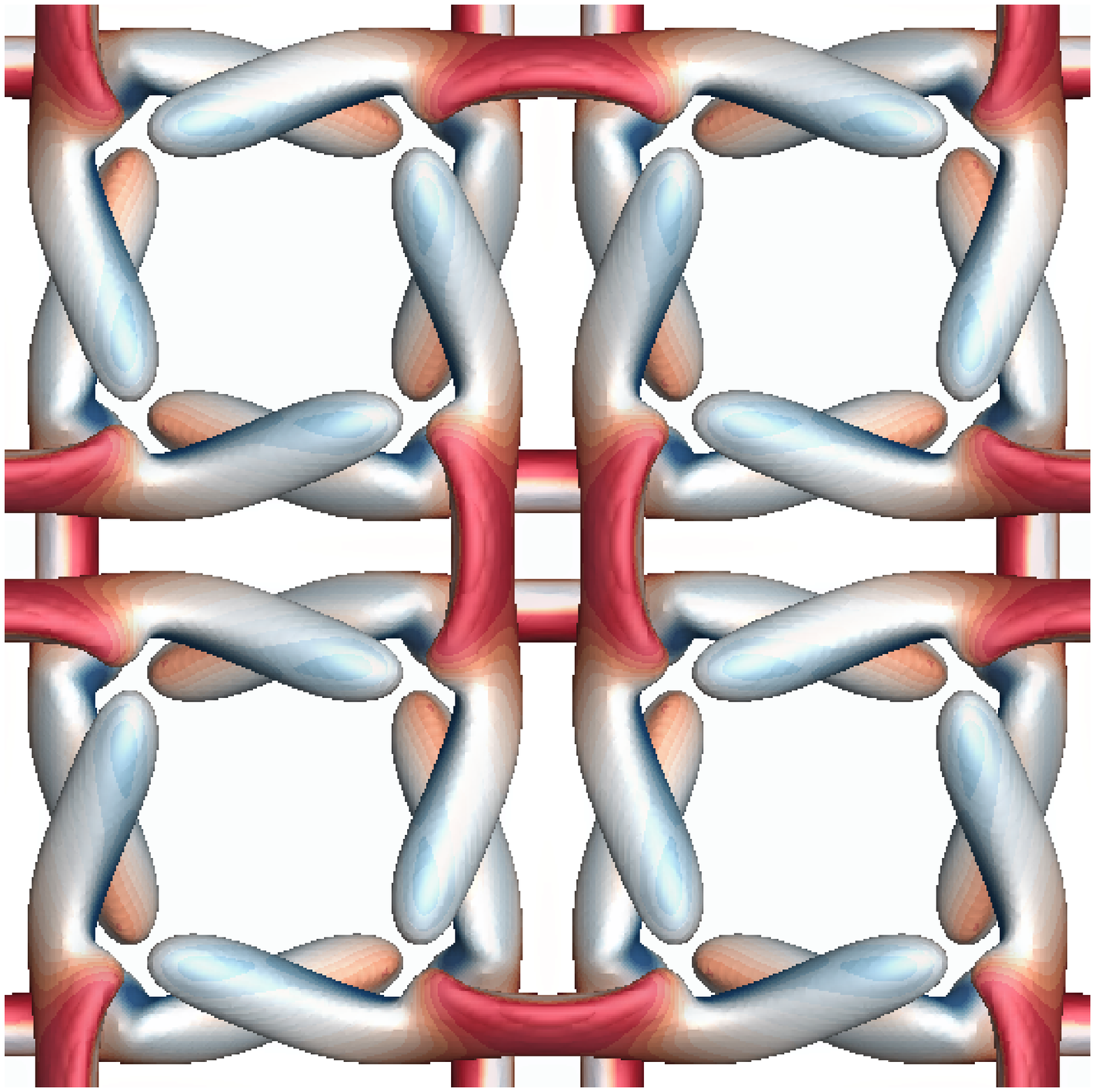}
\includegraphics[width=0.45\textwidth]{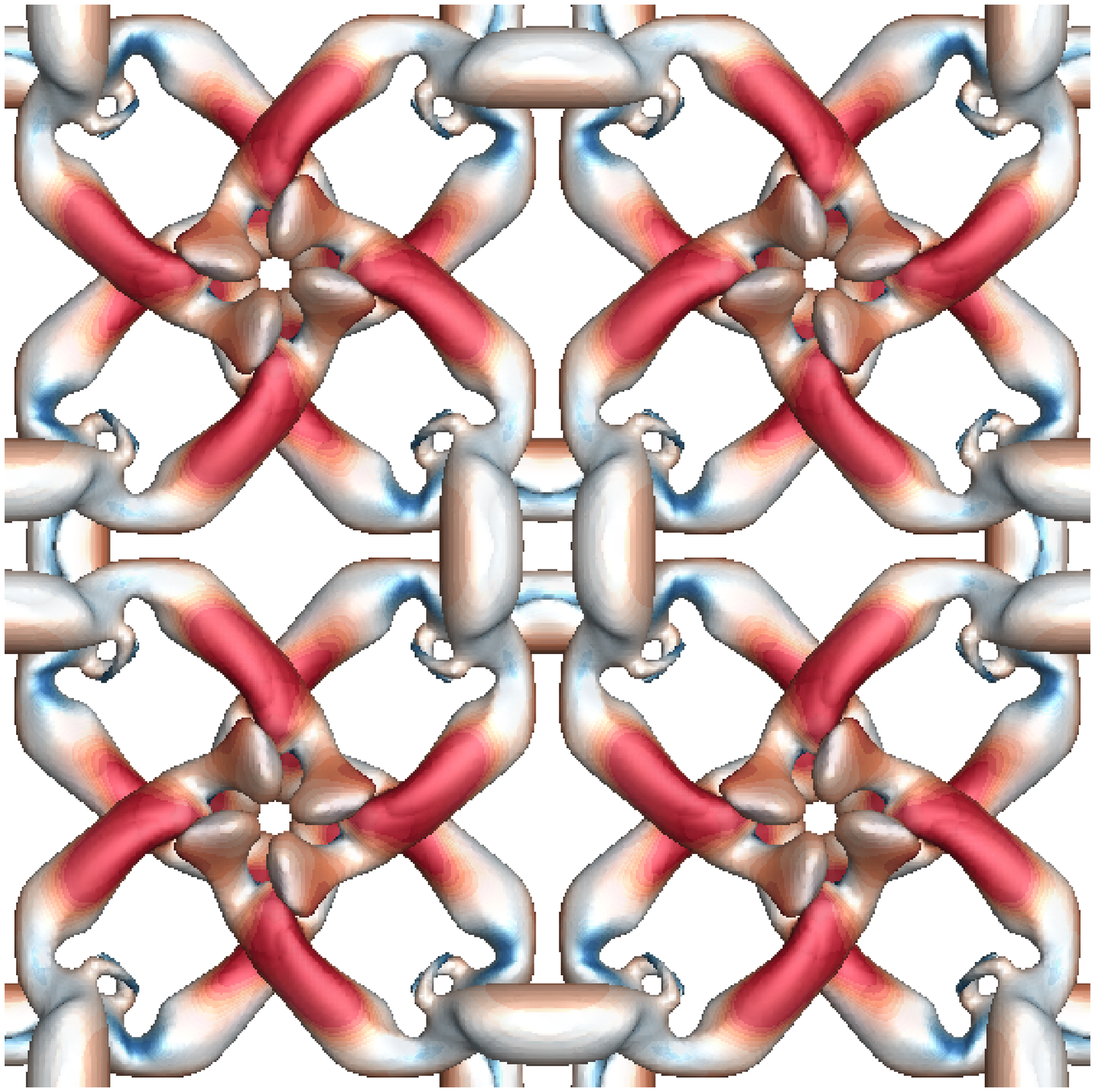}
\caption{\label{tg-vortex-2} Taylor-Green Vortex problem: the
iso-surfaces of $Q$ criterion colored by velocity magnitude at time
$t =5$ and $10$.}
\end{figure}

\subsection{Taylor-Green Vortex}
This problem is aimed at testing the performance of high-order
methods on the direct numerical simulation of a three-dimensional
periodic and transitional flow defined by a simple initial
condition, i.e. the Taylor-Green vortex
\citep{Case-Bull,Case-Debonis}. With a uniform temperature field,
the initial flow field is given by
\begin{align*}
u=&V_0\sin(\frac{x}{L})\cos(\frac{y}{L})\cos(\frac{z}{L}),\\
v=&-V_0\cos(\frac{x}{L})\sin(\frac{y}{L})\cos(\frac{z}{L}),\\
w=&0,\\
p=&p_0+\frac{\rho_0V_0^2}{16}(\cos(\frac{2x}{L})+\cos(\frac{2y}{L}))(\cos(\frac{2z}{L})+2).
\end{align*}
The fluid is then a perfect gas with $\gamma=1.4$ and the Prandtl
number is $Pr=0.71$. Numerical simulations are conducted with two
Reynolds numbers $Re=280$. The flow is computed within a periodic
square box defined as $-\pi L\leq x, y, z\leq \pi L$. The
characteristic convective time $t_c = L/V_0$. In the computation,
$L=1, V_0=1, \rho_0=1$, and the Mach number takes $M_0=V_0/c_0=0.1$,
where $c_0$ is the sound speed. The volume-averaged kinetic energy
can be computed from the flow as it evolves in time, which is
expressed as
\begin{align*}
E_k=\frac{1}{\rho_0\Omega}\int_\Omega\frac{1}{2}\rho\boldsymbol{u}\cdot\boldsymbol{u}\text{d}\Omega,
\end{align*}
where $\Omega$ is the volume of the computational domain, and the
dissipation rate of the kinetic energy is given by
\begin{align*}
\varepsilon_k=-\frac{\text{d}E_k}{\text{d}t}.
\end{align*}
The numerical results with $192\times192\times192$ mesh points for
the normalized volume-averaged kinetic energy and dissipation rate
are presented in Fig.\ref{tg-vortex-1}, which agree well with the
data in \citep{Case-Wang}. The iso-surfaces of $Q$ criterions
colored by velocity magnitude at $t=5$ and $10$ are shown in
Fig.\ref{tg-vortex-2}. The complex structures can be well captured
by the current scheme.

\begin{figure}[!h]
\centering
\includegraphics[width=0.45\textwidth]{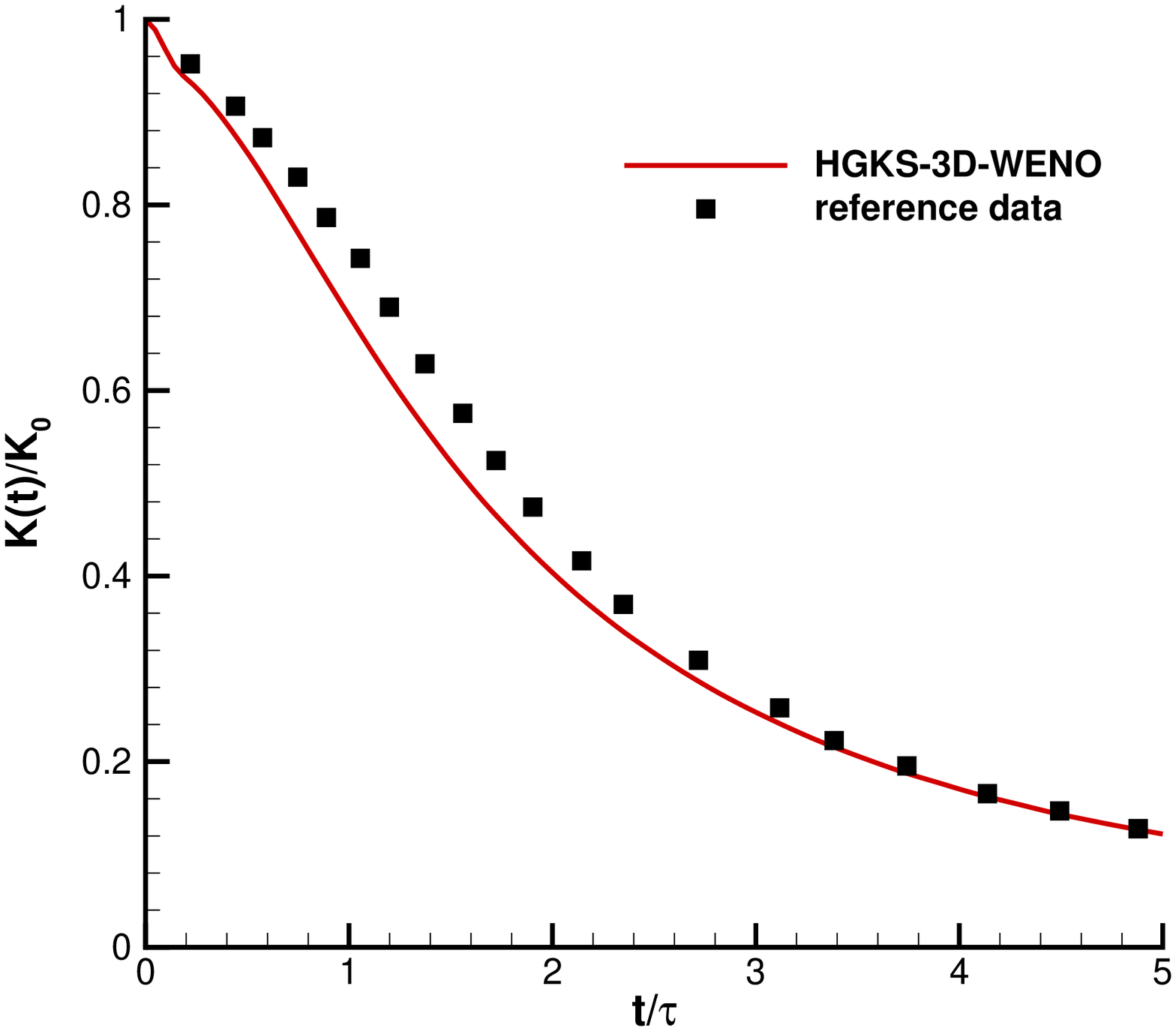}
\includegraphics[width=0.45\textwidth]{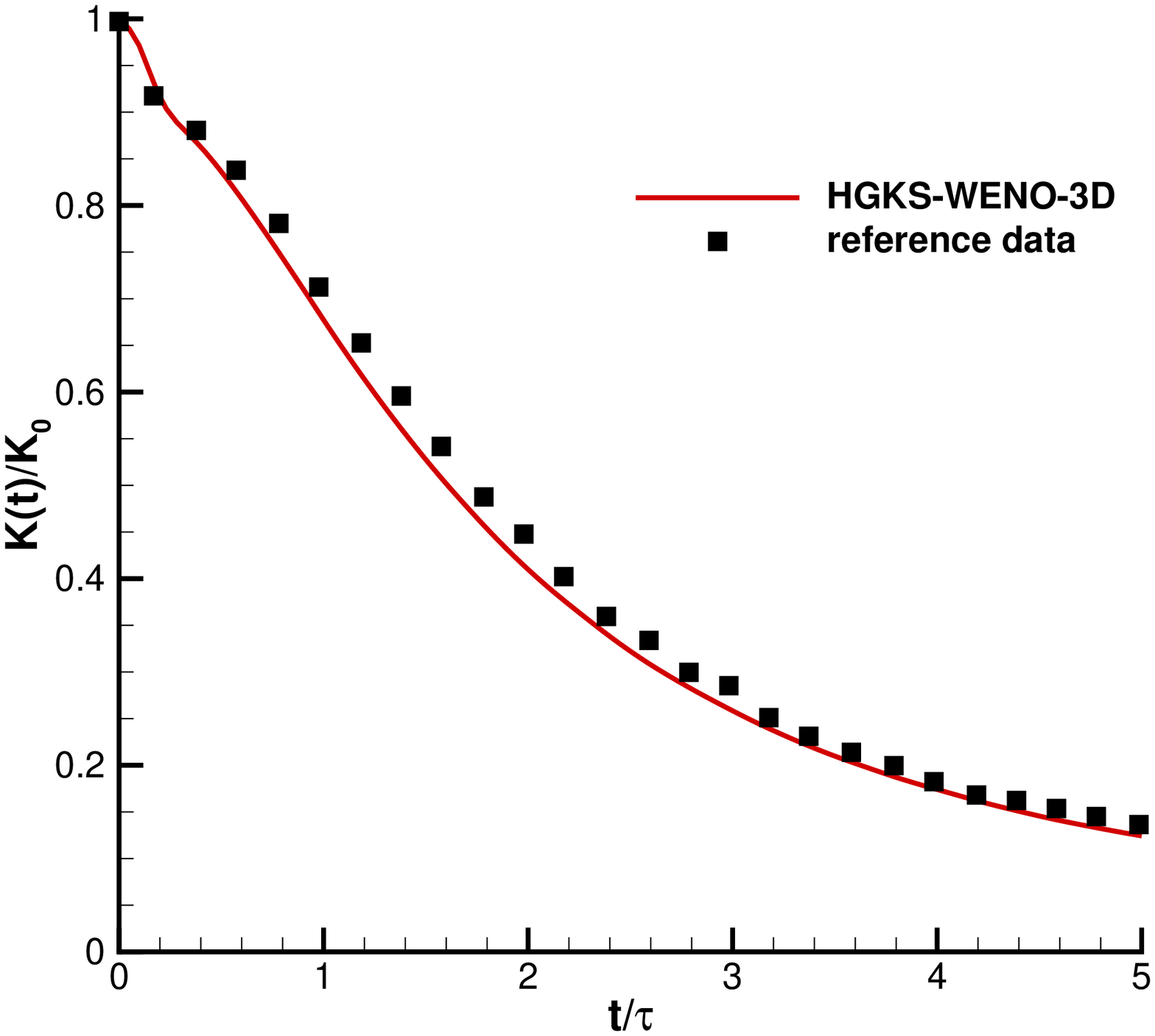}\\
\includegraphics[width=0.45\textwidth]{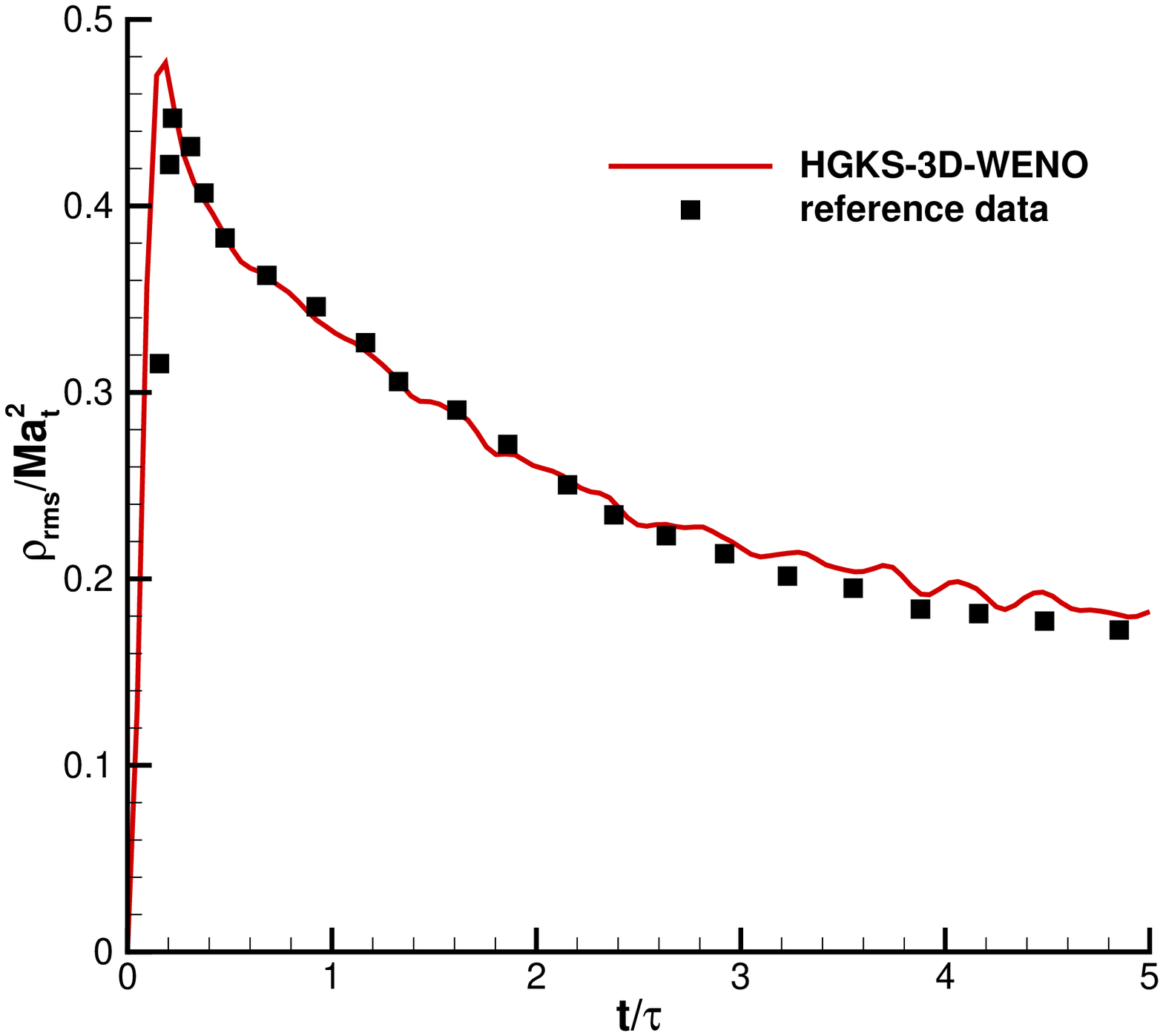}
\includegraphics[width=0.45\textwidth]{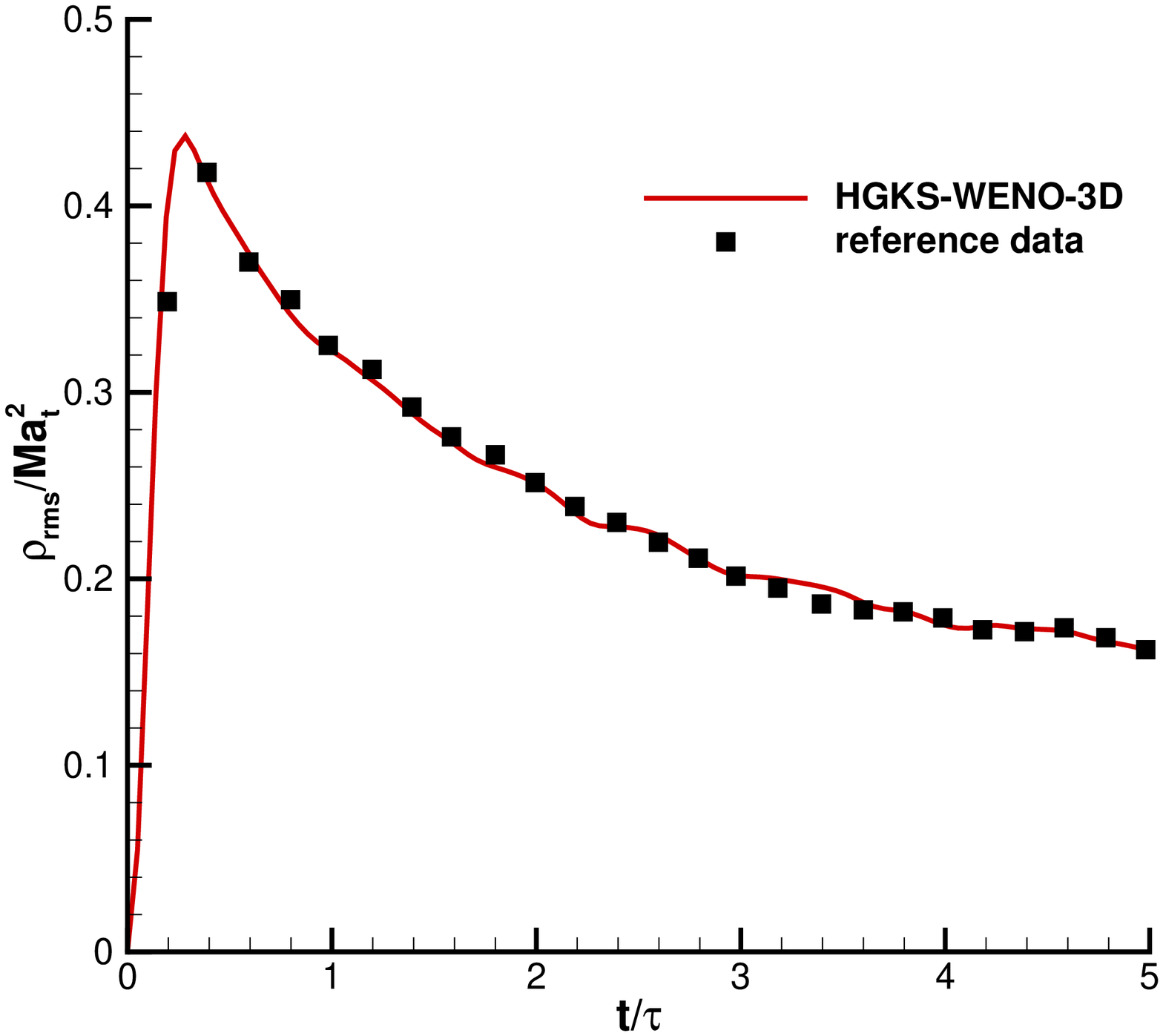}\\
\includegraphics[width=0.45\textwidth]{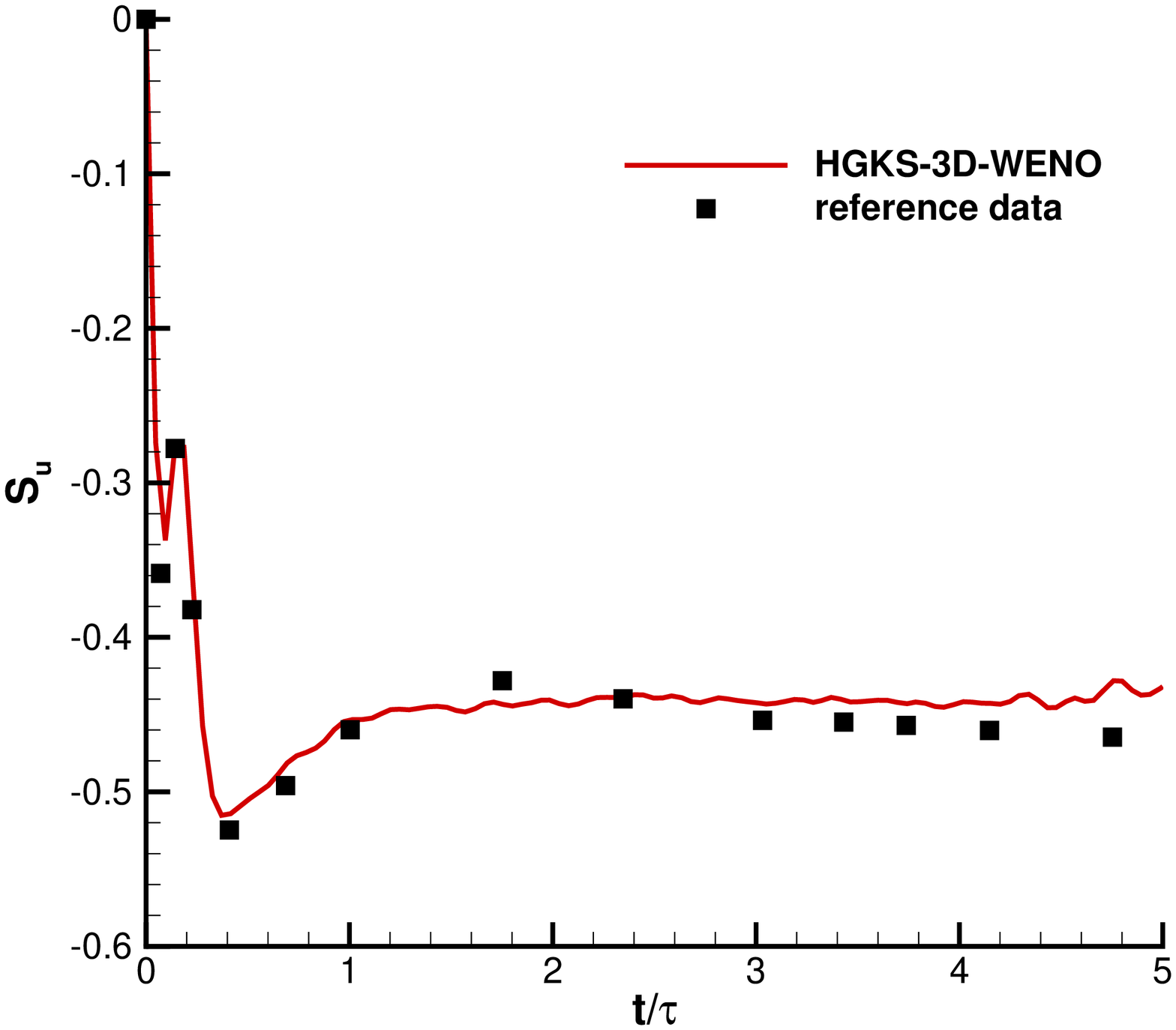}
\includegraphics[width=0.45\textwidth]{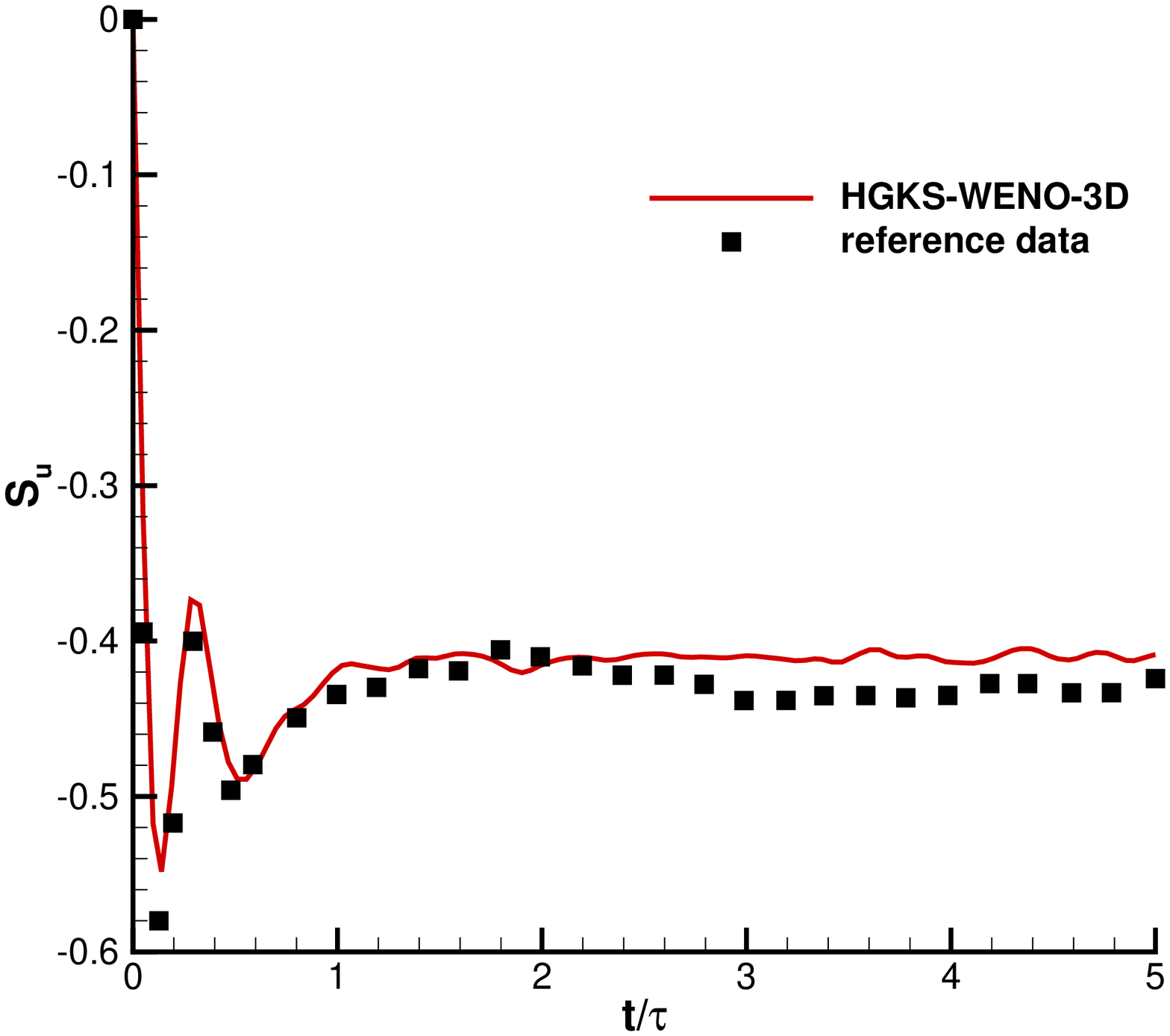}
\caption{\label{homogeneous-1}Compressible isotropic turbulence:
time history of $\rho_{rms}/Ma_t^2$, $K/K_0$, $S_u(t)$ with respect
to $t/\tau$ and Mach number distribution  with $Ma_{t}=0.3$ (left)
and $0.5$ (right).}
\end{figure}

\begin{figure}[!h]
\centering
\includegraphics[width=0.475\textwidth]{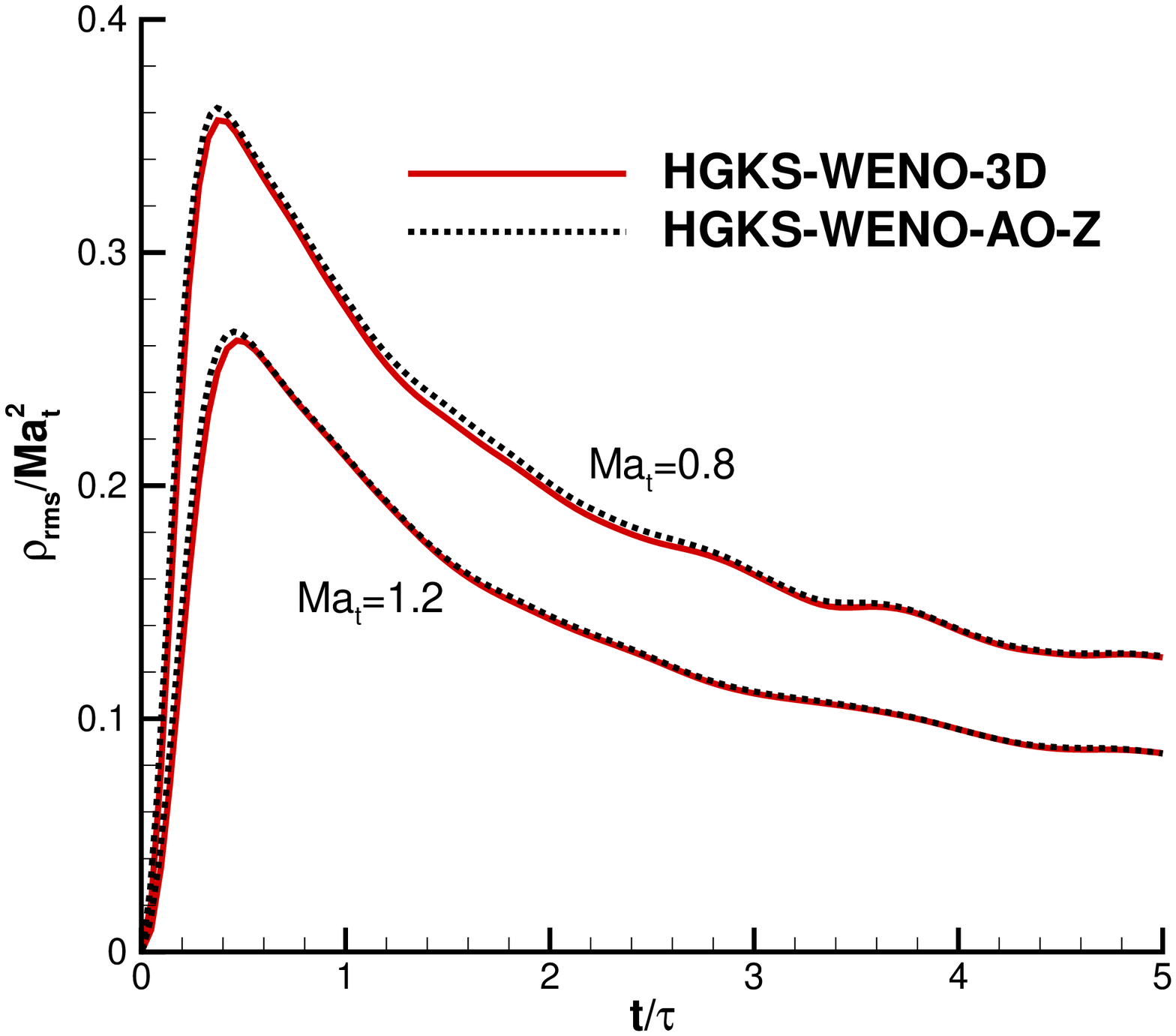}
\includegraphics[width=0.475\textwidth]{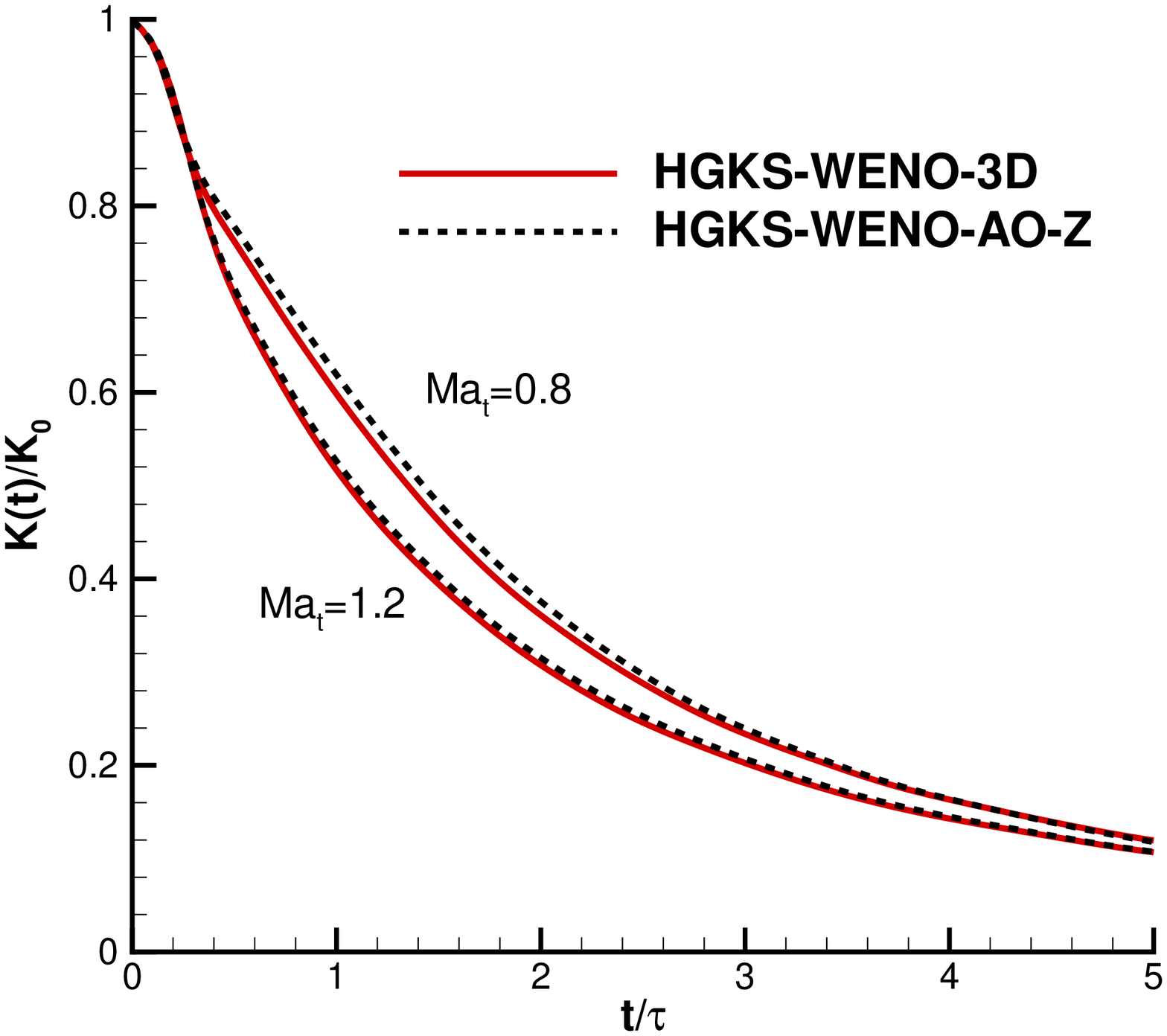}
\caption{\label{homogeneous-2}Compressible isotropic turbulence:
time history of $\rho_{rms}/Ma_t^2$ and $K/K_0$ with respect to
$t/\tau$ with $Ma_{t}=0.8$ and $1.2$.}
\end{figure}

\subsection{Compressible isotropic turbulence}
The compressible isotropic turbulence is regarded as one of
cornerstones to elucidate the effects of compressibility for
compressible turbulence\citep{DNS-2,DNS-3}. The flow is computed
within a square box defined as $-\pi \leq x, y, z\leq \pi$, and the
periodic boundary conditions are used in all directions for all the
flow variables. A divergence-free random initial velocity field
$\boldsymbol{u}_0$ is generated for a given spectrum with a
specified root mean square $u_{rms}$ as follows
\begin{align*}
u_{rms}=<\frac{\boldsymbol{u}\cdot \boldsymbol{u}}{3}>^{1/2},
\end{align*}
where $<...>$ is a volume average over the whole computational
domain. The specified spectrum for velocity is given by
\begin{align*}
E(\kappa)=A_0\kappa^4\exp(-2\kappa^2/\kappa_0^2),
\end{align*}
where $\kappa$ is the wave number, $\kappa_0$ is the wave number at
spectrum peaks, and $A_0$ is a constant chosen to get a specified
initial kinetic energy. With current initial strategy, the initial
ensemble turbulent kinetic energy $K_0$, ensemble enstrophy
$\Omega_0$, ensemble dissipation rate $\varepsilon_0$,
large-eddy-turnover time $\tau_{t_0}$, Kolmogorov length scale
$\eta_0$, and the Kolmogorov time scale $\tau_0$ are given as
\begin{align*}
K_0=&\frac{3A_0}{64} \sqrt{2 \pi} \kappa_0^5, ~ \Omega_0=\frac{15
A_0}{256} \sqrt{2 \pi} \kappa_0^7,~
\tau_{t_0}=\sqrt{\frac{32}{A_0}}(2 \pi)^{1/4} \kappa_0^{-7/2},\\
&\varepsilon_0=\frac{2\mu_0\Omega_0}{\rho_0}, ~
\eta_0=(\nu_0^3/\varepsilon_0)^{1/4}, ~
\tau_0=(\nu_0/\varepsilon_0)^{1/2}.
\end{align*}
The evolution of this system is dominated by the initial
thermodynamic quantities and two dimensionless parameters, i.e. the
initial Taylor microscale Reynolds number and turbulent Mach number
\begin{align*}
Re_\lambda=&\frac{<\rho>u_{rms}\lambda}{<\mu>}=\frac{(2\pi)^{1/4}}{4}\frac{\rho_0}{\mu_0}\sqrt{2A_0}k_0^{3/2},\\
&Ma_t=\frac{\sqrt{3}u_{rms}}{<c_s>}=\frac{\sqrt{3}u_{rms}}{\sqrt{\gamma
T_0}},
\end{align*}
where $\lambda$ is Taylor microscale
\begin{align*}
\lambda^2=\frac{u_{rms}^2}{<(\partial_x u)^2>}.
\end{align*}
The dynamic viscosity is determined by
\begin{align*}
\mu=\mu_0\big(\frac{T}{T_0}\big)^{0.76},
\end{align*}
where $\mu_0$ and $T_0$ can be determined from $Re_\lambda$ and
$Ma_t$ with initialized $u_{rms}$ and $\rho_0=1$.

\begin{figure}[!h]
\centering
\includegraphics[width=0.475\textwidth]{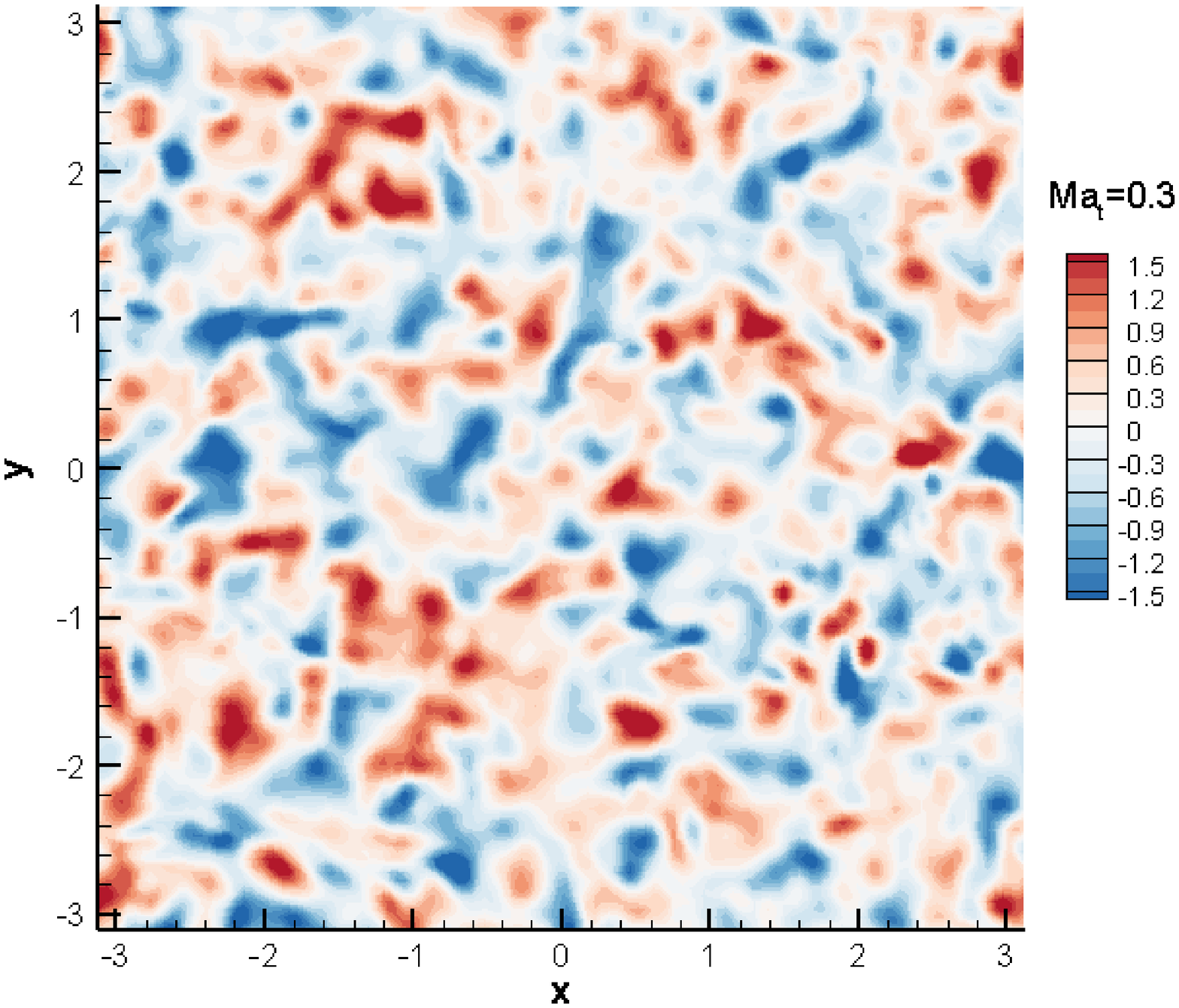}
\includegraphics[width=0.475\textwidth]{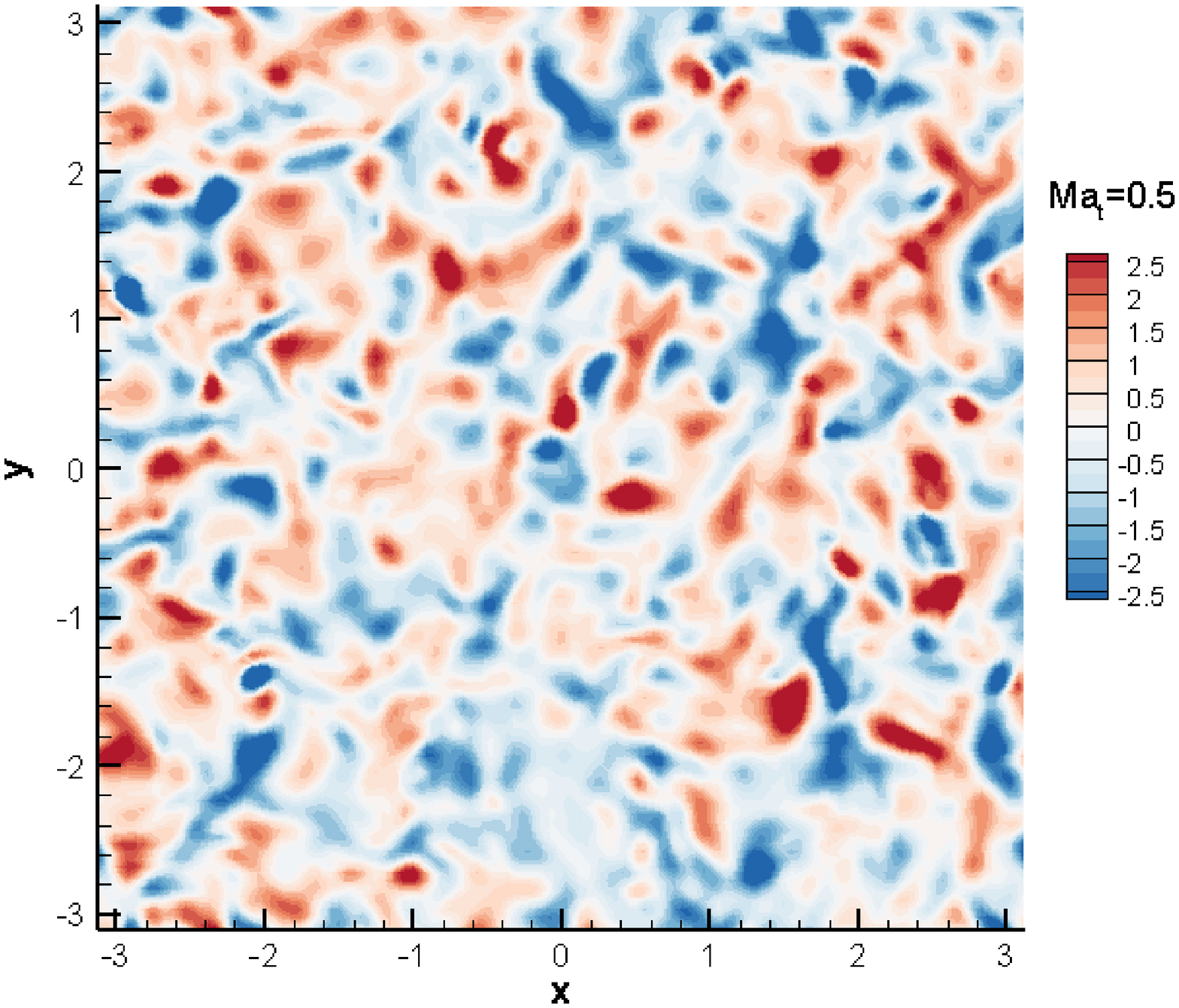}
\includegraphics[width=0.475\textwidth]{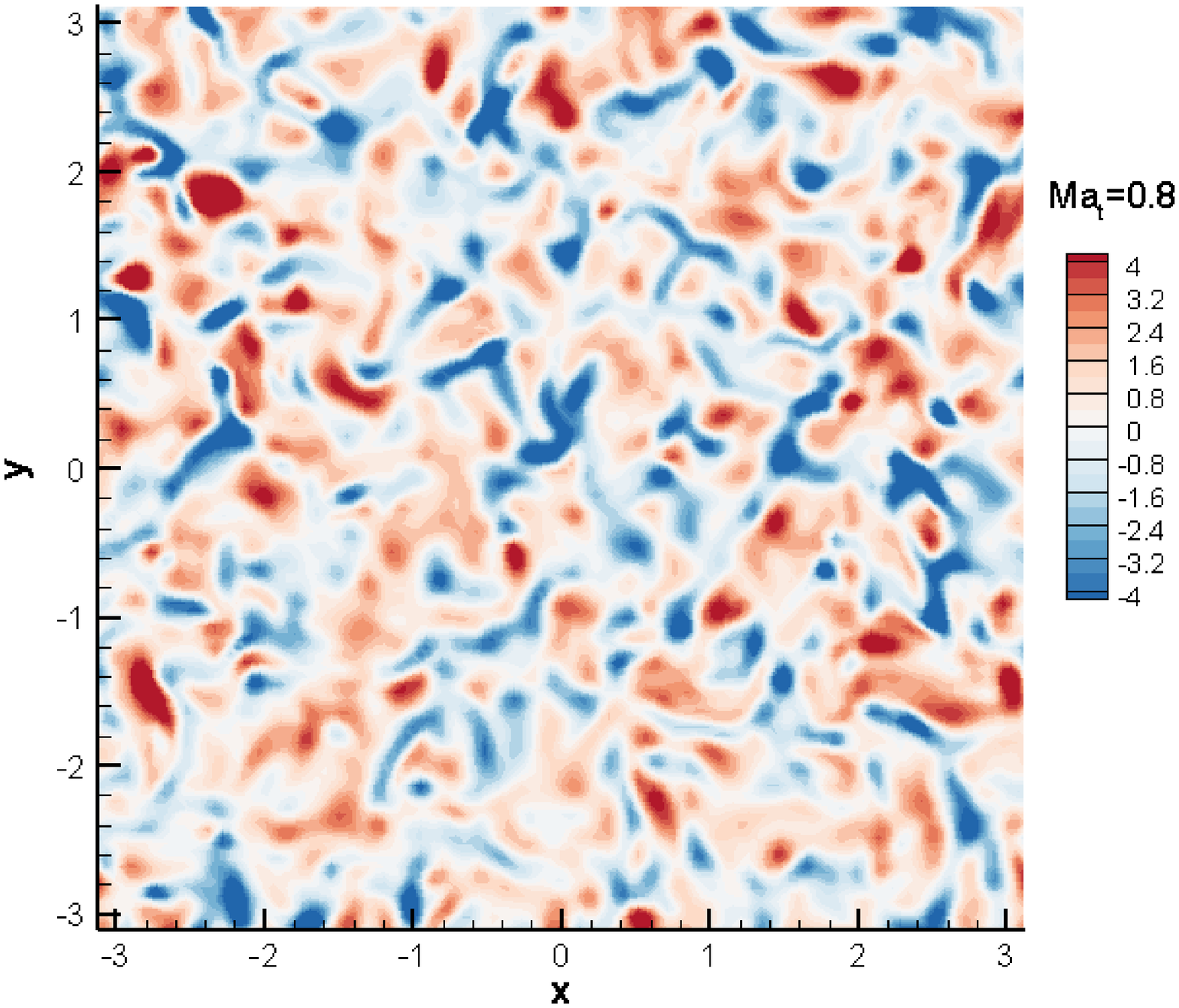}
\includegraphics[width=0.475\textwidth]{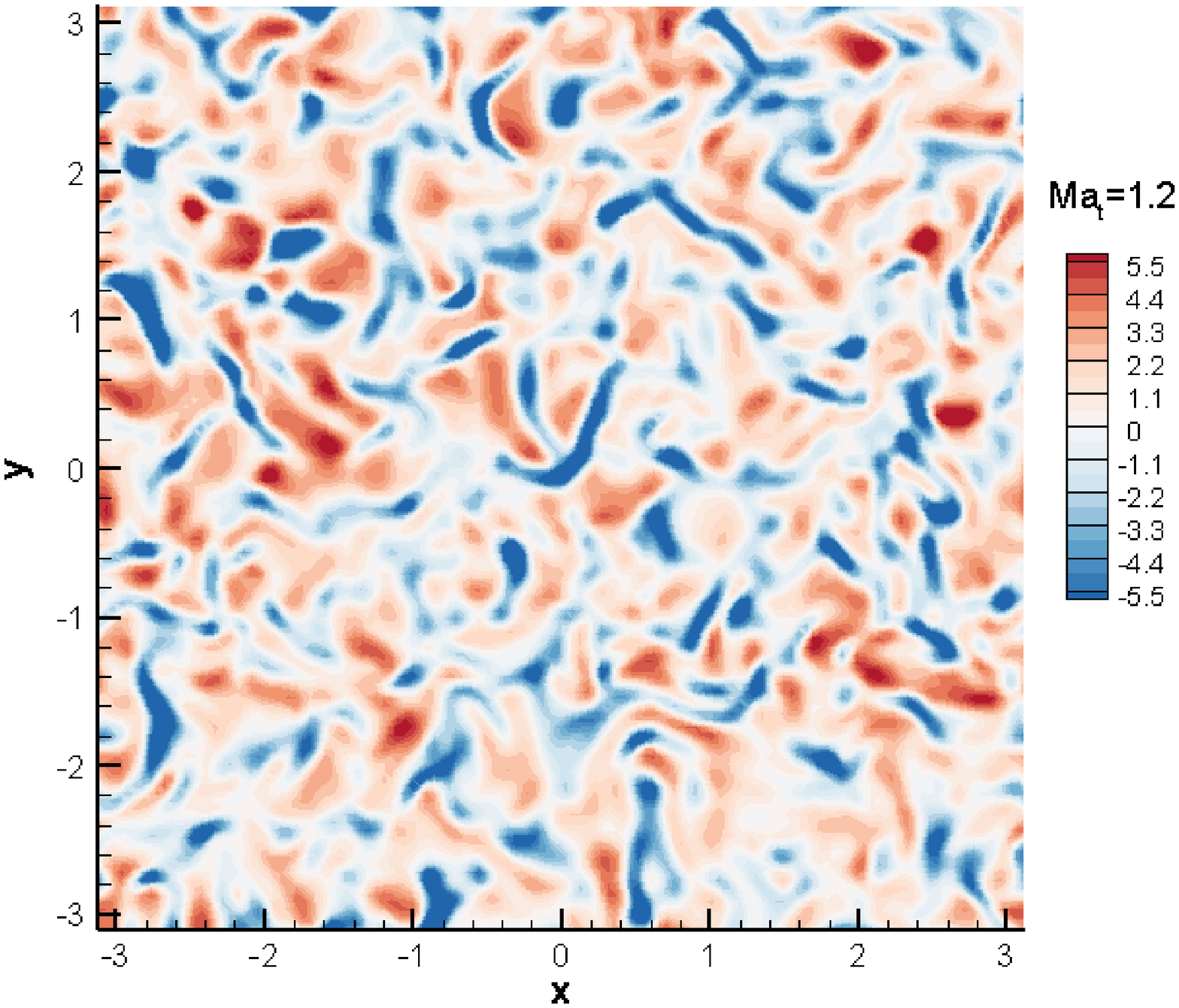}
\caption{\label{homogeneous-3}Compressible isotropic turbulence: the
contours of dilation $\theta=\nabla\cdot \boldsymbol{u}$ at
$t/\tau=1$ with $Ma_{t}=0.3, 0.5, 0.8$ and $1.2$.}
\end{figure}

High-order compact finite difference method \cite{DNS-1} has been
widely utilized in the simulation of isotropic compressible
turbulence with moderate turbulent Mach number ($Ma_t \leq 0.8$).
However, when simulating the turbulent in supersonic regime ($Ma_t
\geq 0.8$), the above scheme fails to capture strong shocklets and suffers from
numerical instability. However, the current HGKS scheme will be tested at a wide range
turbulent Mach numbers. In the computation,
$A_0=1.3\times10^{-4}, k_0=8$, $Re_\lambda=72$ and the uniform
meshes with $128^3$ cells are used. The compressible isotropic
turbulent flows in nonlinear subsonic regime with $Ma_{t}=0.3$ and
$0.5$ are tested firstly. The time history of normalized kinetic
energy $K(t)/K_0$, normalized root-mean-square of density
fluctuation $\rho_{rms}(t)/Ma_t^2$ and skewness factor $S_u(t)$ with
respect to $t/\tau$ are given in Fig.\ref{homogeneous-1}.  The
numerical results agree well with the reference data \citep{DNS-2}.
With fixed initial $Re_\lambda=72$ and $128^3$ cells, the cases with
$Ma_t=0.8$ and $1.2$ are tested as well, which go up to the
supersonic turbulent Mach number. The time histories of normalized
kinetic energy $K(t)/K_0$ and normalized root-mean-square of density
fluctuation $\rho_{rms}(t)/Ma_t^2$ at $t/\tau=1$ are given in
Fig.\ref{homogeneous-2} as well. With the increase of $Ma_t$, the
dynamic viscosity increases and the kinetic energy gets dissipated
more rapidly. As comparison, the numerical results with fifth-order
WENO-Z scheme are given as well. More studies of compressible
isotropic turbulence can be referred in \cite{GKS-high-4}. The
contours of dilation $\theta=\nabla\cdot \boldsymbol{u}$ for
$Ma_t=0.3, 0.5, 0.8$ and $1.2$ are given in Fig.\ref{homogeneous-3},
which shows very different behavior between the compression motion
and expansion motion. With the increase of $Ma_t$, the compression
regions, i.e. shocklets behave in the shape of narrow and long
``ribbon". In addition, the strong compression regions are close to
several regions of high expansion.  Compared the case with
$Ma_t=0.3$ in subsonic regime, the supersonic case with $Ma_t=1.2$
contains much more crisp shocklets, which pose much greater
challenge for high-order schemes when implementing DNS for isotropic
turbulence in supersonic regime, which validate the robustness for
the challenging compressible turbulence problems.

\section{Conclusion}
In this paper, with the WENO-AO reconstruction an efficient and simple third-order gas-kinetic scheme  is developed
for the three-dimensional Euler and Navier-Stokes equations. In the classical WENO
scheme, choosing sub-stencils from big stencil and solving linear
weights at Gaussian quadrature points would make the reconstruction
complicated, especially for three-dimensional flows. To overcome the
drawback, the WENO-AO strategy is adopted. Based on the candidate
stencils, the quadratic polynomial for big stencil and linear
polynomials for sub-stencils are constructed.  The spatial
independent linear weights are used, which have fixed values and
become positive. With the smooth indicator, the nonlinear weights
can be constructed. Through particle colliding procedure, the
point-value and slopes for equilibrium part are obtained directly from the initial
reconstruction of the non-equilibrium state, an extra reconstruction for the equilibrium state in
the classical HGKS is avoided. Taken the grid velocity into account,
such scheme can be also extended to the moving-mesh computation.
For the mesh with non-coplanar vertexes, which is commonly generated
in the moving-mesh computation, the trilinear interpolation is used
to parameterize the hexahedron, and the bilinear interpolation is
used to parameterize the interface of hexahedron. Numerical results
are provided to illustrate the good performance of the WENO schemes
from the smooth inviscid flows to the supersonic turbulent flows. In
the future, the extension to unstructured meshes will be developed.

\section*{Acknowledgements}
The current research of L. Pan is supported by National Science
Foundation of China (11701038) and the Fundamental Research Funds
for the Central Universities. The work of K. Xu is supported by
National Science Foundation of China (11772281, 91852114) and Hong
Kong research grant council (16206617).

\end{document}